\documentclass[10pt]{amsart}

\usepackage{url,graphicx}
\usepackage[dvips]{epsfig}
\usepackage{latexsym,color}
\usepackage{amscd,amssymb,verbatim}

\newcommand{\bbr}{\big{(}\big{)}}

\newcommand{\be}{\begin{enumerate}}

\newcommand{\bo}{\partial}
\newcommand{\bu}{\bullet}

\newcommand{\cc}{{\mathcal C}}

\newcommand{\chom}{\text{\tt Hom}_{\,0}}

\newcommand{\cn}{{\mathcal N}}

\newcommand{\com}{\complement}

\newcommand{\cs}{{\mathcal S}}

\newcommand{\da}{\Delta}

\newcommand{\dc}{{\mathbb C}}
\newcommand{\dr}{{\mathbb R}}

\newcommand{\dz}{{\mathbb Z}}

\newcommand{\ee}{\end{enumerate}}

\newcommand{\im}{{\text{\rm Im}}\,}
\newcommand{\ind}{{\text{\rm Ind}}\,}

\newcommand{\lra}{\longrightarrow}

\newcommand{\nin}{\noindent}

\newcommand{\pr}{\noindent{\bf Proof. }}

\newcommand{\ra}{\rightarrow}

\newcommand{\rp}{{\mathbb R}{\mathbb P}}

\newcommand{\sgn}{\text{sgn}\,}

\newcommand{\sm}{\setminus}

\newcommand{\supp}{\text{\rm supp}\,}
\newcommand{\susp}{\text{\rm susp}}

\newcommand{\thom}{\text{\tt Hom}\,}
\newcommand{\thomp}{\text{\tt Hom}_+}
\newcommand{\ti}{\tilde }

\newcommand{\vt}{\vartheta}

\newcommand{\wti}{\widetilde }

\newcommand{\zz}{{{\mathbb Z}_2}}

\newtheorem{thm}{Theorem}[section]
\newtheorem{df}  [thm]{Definition}

\newtheorem{lm}  [thm]{Lemma}
\newtheorem{crl} [thm]{Corollary}
\newtheorem{prop}[thm]{Proposition}
\newtheorem{conj}[thm]{Conjecture}
\newtheorem{rem} [thm]{Remark}

\numberwithin{equation}{section}
\numberwithin{figure}{section}
\numberwithin{table}{section}

\begin{document}

\title {Proof of the Lov\'asz Conjecture}

\author{Eric Babson and Dmitry N. Kozlov}

\date{\noindent
  \today\\[0.05cm]
  }

\address{Department of Mathematics, University of Washington, Seattle,
  U.S.A.} \email{babson@math.washington.edu}

\address{Institute of Theoretical Computer Science, Eidgen\"ossische 
Technische Hochschule, Z\"urich, Switzerland}
\email{dkozlov@inf.ethz.ch}

\begin{abstract}
To any two graphs $G$ and $H$ one can associate a~cell complex
$\thom(G,H)$ by taking all graph multihomorphisms from $G$ to $H$ as
cells.

In this paper we prove the Lov\'asz Conjecture which states that
\begin{quote}
{\it if $\thom(C_{2r+1},G)$ is $k$-connected, then $\chi(G)\geq k+4$,}
\end{quote}
where $r,k\in\dz$, $r\geq 1$, $k\geq -1$, and $C_{2r+1}$ denotes the
cycle with $2r+1$ vertices.

The proof requires analysis of the complexes $\thom(C_{2r+1},K_n)$.
For even $n$, the obstructions to graph colorings are provided by the
presence of torsion in $H^*(\thom(C_{2r+1},K_n);\dz)$. For odd $n$,
the obstructions are expressed as vanishing of certain powers of
Stiefel-Whitney characteristic classes of $\thom(C_{2r+1},K_n)$, where
the latter are viewed as $\zz$-spaces with the involution induced by
the reflection of $C_{2r+1}$.
\end{abstract}

\maketitle

\section{Introduction}\label{sect_intr}

The main idea of this paper is to look for obstructions to graph
colorings in the following indirect way: take a~graph, associate to it
a~topological space, and then look for obstructions to colorings of
the graph by studying the algebraic invariants of this space.

The construction of such a space, which is of interest here, has been
suggested by L.\ Lov\'asz.  The obtained complex $\thom(G,H)$ depends
on two graph parameters. The algebraic invariants of this space, which
we proceed to study, are its cohomology groups, and, when it can be
viewed as a~$\zz$-space, its Stiefel-Whitney characteristic classes.

\subsection{The vertex colorings and the category of Graphs.} 
$\,$ \vspace{5pt}

All graphs in this paper are undirected. The following definition is
a~key in turning the set of all undirected graphs into a~category.

\begin{df}\label{dfgrhom}
  For two graphs $G$ and $H$, a {\bf graph homomorphism}
  from $G$ to $H$ is a~map $\phi:V(G)\rightarrow V(H)$, such that if
  $(x,y)\in E(G)$, then $(\phi(x),\phi(y))\in E(H)$.
\end{df}
Here, $V(G)$ denotes the set of vertices of $G$, and $E(G)$ denotes
the set of its edges.

For a graph $G$ the {\it vertex coloring} is an assignment of colors
to vertices such that no two vertices which are connected by an~edge
get the same color. The minimal needed number of colors is denoted by
$\chi(G)$, and is called the {\it chromatic number} of~$G$.

Deciding whether or not there exists a~graph homomorphism between two
graphs is in general at least as difficult as bounding the chromatic
numbers of graphs because of the following observation: a vertex
coloring of $G$ with $n$ colors is the same as a graph homomorphism
from $G$ to the complete graph on $n$ vertices $K_n$. Because of this,
one can also think of graph homomorphisms from $G$ to $H$ as vertex
colorings of $G$ with colors from $V(H)$ subject to the natural
condition.

Since an identity map is a~graph homomorphism, and a~composition of
two graph homomorphisms is again a~graph homomorphism, we can consider
the category {\bf Graphs} whose objects are all undirected graphs, and
morphisms are all the graph homomorphisms.

We denote the set of all graph homomorphisms from $G$ to $H$ by
$\chom(G,H)$. Lov\'asz has suggested the following way of turning this
set into a~topological space.

\begin{df} \label{dfhom}
We define $\thom(G,H)$ to be a~polyhedral complex whose cells are
indexed by all functions $\eta:V(G)\rightarrow
2^{V(H)}\setminus\{\emptyset\}$, such that if $(x,y)\in E(G)$, for any
$\tilde x\in\eta(x)$ and $\tilde y\in\eta(y)$ we have $(\tilde
x,\tilde y)\in E(H)$.

  The closure of a~cell $\eta$ consists of all cells indexed by
  $\ti\eta:V(G)\rightarrow 2^{V(H)}\setminus\{\emptyset\}$, which
  satisfy $\ti\eta(v)\subseteq\eta(v)$, for all $v\in V(G)$.
\end{df}
\nin We think of a cell in $\thom(G,H)$ as a~collection of non-empty
lists of vertices of $H$, one for each vertex of $G$, with the
condition that any choice of one vertex from each list will yield
a~graph homomorphism from $G$ to~$H$. A~geometric realization of
$\thom(G,H)$ can be described as follows: number the vertices of $G$
with $1,\dots,|V(G)|$, the cell indexed with $\eta:V(G)\rightarrow
2^{V(H)}\setminus\{\emptyset\}$ is realized as a~direct product of
simplices $\da^1,\dots,\da^{|V(G)|}$, where $\da^i$ has $|\eta(i)|$
vertices and is realized as the standard simplex in ${\mathbb
  R}^{|\eta(i)|}$. In particular, the set of vertices of $\thom(G,H)$
is precisely $\chom(G,H)$.

The barycentric subdivision of $\thom(G,H)$ is isomorphic as
a~simplicial complex to the geometric realization of its face poset.
So, alternatively, it could be described by first defining a~poset of
all $\eta$ satisfying conditions of Definition~\ref{dfhom}, with
$\eta\geq\ti\eta$ iff $\eta(v)\supseteq\ti\eta(v)$, for all $v\in
V(G)$, and then taking the geometric realization.

 The $\thom$ complexes are functorial in the following sense:
 $\thom(H,-)$ is a covariant, while $\thom(-,H)$ is
a~contravariant functor from {\bf Graphs} to {\bf Top}. If
$\phi\in\chom(G,G')$, then we shall denote the induced cellular
maps as $\phi^H:\thom(H,G)\ra\thom(H,G')$ and
$\phi_H:\thom(G',H)\ra\thom(G,H)$.

\subsection{The statement of the Lov\'asz conjecture.} $\,$ \vspace{5pt}

Lov\'asz has stated the~following conjecture, which we prove in this
paper.

\begin{thm} {\bf (Lov\'asz Conjecture).} \label{loconj}
  {\it Let $G$ be a graph, such that $\thom(C_{2r+1},G)$ is
    $k$-connected for some $r,k\in\dz$, $r\geq 1$, $k\geq -1$, then
    $\chi(G)\geq k+4$.}
\end{thm}

\nin Here $C_{2r+1}$ is a cycle with $2r+1$ vertices:
$V(C_{2r+1})=\dz_{2r+1}$,
$E(C_{2r+1})=\{(x,x+1),(x+1,x)\,|\,x\in\dz_{2r+1}\}$. 

The motivation for this conjecture stems from the following theorem
which Lov\'asz has proved in 1978.

\begin{thm} \label{lothm}
  {\rm (Lov\'asz, \cite{Lo}).}  {\it Let $H$ be a graph, such that
    $\thom(K_2,H)$ is $k$-connected for some $k\in\dz$, $k\geq -1$,
    then $\chi(H)\geq k+3$.}
\end{thm}

One corollary of Theorem~\ref{lothm} is the Kneser conjecture from
1955, see~\cite{Knes}.

\begin{rem} The actual theorem from \cite{Lo} is stated using the
neighborhood complexes $\cn(H)$. However, it is well known that
$\cn(H)$ is homotopy equivalent to $\thom(K_2,H)$ for any graph~$H$,
see, e.g., \cite{BK2} for an argument. In fact, these two spaces are
known to be simple-homotopy equivalent, see \cite{K5}.
\end{rem}

We note here that Theorem~\ref{loconj} is trivially true for $k=-1$:
$\thom(C_{2r+1},G)$ is $(-1)$-connected if and only if it is
non-empty, and since there are no homomorphisms from odd cycles to
bipartite graphs, we conclude that $\chi(G)\geq 3$. It is also not
difficult to show that Theorem~\ref{loconj} holds for $k=0$ by using
the winding number. A short argument for a~more general statement can
be found in subsection~\ref{ss2.2}.

\subsection{Plan of the paper.} $\,$
\vspace{5pt}

In Section 2, we formulate the main theorems and describe the general
framework of finding obstructions to graph colorings via vanishing of
powers of Stiefel-Whitney characteristic classes.

In Section 3, we introduce auxiliary simplicial complexes, which we
call $\thomp(-,-)$. For any two graphs $G$ and $H$, there is
a~canonical support map $\supp:\thomp(G,H)\ra\Delta_{|V(G)|-1}$, and
the preimage of the barycenter is precisely $\thom(G,H)$.  This allows
us to set up a~useful spectral sequence, filtering by the preimages of
the $i$-skeleta.

In Section 4, we compute the cohomology groups
$H^*(\thom(C_{2r+1},K_n);\dz)$ up to dimension $n-2$, and we find the
$\zz$-action on these groups. These computations allow us to prove the
Lov\'asz conjecture for the case of odd~$k$, $k\geq 1$.

In Section 5, we study a~different spectral sequence, this one
converging to $H^*(\thom(C_{2r+1},K_n)/\zz;\zz)$. Understanding
certain entries and differentials leads to the proof of the Lov\'asz
conjecture for the case of even~$k$ as well.

The results of this paper were announced in \cite{BK1}, where no
complete proofs were given. The reader is referred to \cite{IAS} for
a~survey on $\thom$ complexes, which also includes a lot of background
material which is omitted in this paper.
\vspace{5pt}

\nin {\bf Acknowledgments.} The second author acknowledges support by
the University of Washington, Seattle, the Swiss National Science
Foundation, and the Swedish National Research Council.

\section{The idea of the proof of the Lov\'asz conjecture.}

\subsection{Group actions on $\thom$ complexes and Stiefel-Whitney 
classes.}
$\,$  \vspace{5pt} \label{ss2.1}

Consider an arbitrary CW complex $X$ on which a~finite group $\Gamma$
acts freely. By the general theory of principal $\Gamma$-bundles,
there exists a~$\Gamma$-equivariant map $\tilde w:X\ra{\bf E}\Gamma$,
and the induced map $w:X/\Gamma\ra{\bf B}\Gamma={\bf E}\Gamma/\Gamma$
is unique up to homotopy.

Specifying $\Gamma=\zz$, we get a map $\tilde w:X\ra S^{\infty}={\bf
E}\zz$, where $\zz$ acts on $S^{\infty}$ by the antipodal map, and the
induced map $w:X/\zz\ra{\mathbb R\mathbb P}^\infty={\bf B}\zz$.  We
denote the induced $\zz$-algebra homomorphism $H^*({\mathbb R \mathbb
P}^\infty;\zz)\ra H^*(X/\zz;\zz)$ by $w^*$. Let $z$ denote the
nontrivial cohomology class in $H^1({\mathbb R\mathbb P}^\infty;
\zz)$. Then $H^*({\mathbb R\mathbb P}^\infty;\zz) \simeq\zz[z]$ as
a~graded $\zz$-algebra, with $z$ having degree~1. We denote the image
$w^*(z)\in H^1(X/\zz;\zz)$ by $\varpi_1(X)$. This is the {\it first
Stiefel-Whitney class} of the $\zz$-space~X. Clearly,
$\varpi_1^k(X)=w^*(z^k)$, since $w^*$ is a~$\zz$-algebra homomorphism.
We will be mainly interested in the {\it height} of the
Stiefel-Whitney class, i.e., largest $k$, such that $\varpi_1^k(X)\neq
0$; it was called cohomology co-index in~\cite{CF}.

Turning to graphs, let $G$ be a graph with $\zz$-action given by
$\phi:G\ra G$, $\phi\in\chom(G,G)$, such that $\phi$ flips an edge,
that is, there exist $a,b\in V(G)$, $a\neq b$, $(a,b)\in E(G)$, such
that $\phi(a)=b$ (which implies $\phi(b)=a$). For any graph $H$ we
have the induced $\zz$-action $\phi_H:\thom(G,H)\ra\thom(G,H)$. In
case $H$ has no loops, it follows from the fact that $\phi$ flips an
edge that this $\zz$-action is free.

Indeed, since $\phi_H$ is a~cellular map, if it fixes a~point from
some cell $\eta:V(G)\ra 2^{V(H)}\sm\{\emptyset\}$, then it maps $\eta$
onto itself. By definition, $\phi$ maps $\eta$ to $\eta\circ\phi$, so
this means that $\eta=\eta\circ\phi$. In particular,
$\eta(a)=\eta\circ\phi(a)=\eta(b)$. Since $\eta(a)\neq\emptyset$, we
can take $v\in V(H)$, such that $v\in\eta(a)$. Now, $(a,b)\in E(G)$,
but $(v,v)\notin E(H)$, since $H$ has no loops, which contradicts the
fact that $\eta\in\thom(G,H)$.

Therefore, in this situation, $\thom(G,-)$ is a~covariant functor from
the induced subcategory of {\bf Graphs}, consisting of all loopfree
graphs, to {\bf $\zz$-spaces} (the category whose objects are
$\zz$-spaces and morphisms are $\zz$-maps).

We order $V(C_{2r+1})$ by identifying it with $[1,2r+1]$ by the map
$q:\dz\ra\dz_{2r+1}$, taking $x\mapsto[x]_{2r+1}$. With this notation
$\zz$ acts on $C_{2r+1}$ by mapping $[x]_{2r+1}$ to $[-x]_{2r+1}$, for
$x\in V(C_{2r+1})$. Let $\gamma\in\chom(C_{2r+1},C_{2r+1})$ denote the
corresponding graph homomorphism. This action has a~fixed point
$2r+1$, and it flips one edge $(r,r+1)$.

Furthermore, let $\zz$ act on $K_m$ for $m\geq 2$, by swapping the
vertices 1 and 2 and fixing the vertices $3,\dots,m$; here, $K_m$ is
the graph defined by $V(K_m)=[1,m]$, $E(K_m)=\{(x,y)\,|\,x,y\in V(K_m),
x\neq y\}$. Since in both cases the graph homomorphism flips an~edge,
they induce free $\zz$-actions on $\thom(C_{2r+1},G)$ and
$\thom(K_m,G)$, for an arbitrary graph $G$ without loops.

\subsection{Non-vanishing of powers of Stiefel-Whitney classes as
  obstructions to graph colorings.} $\,$ \vspace{5pt}
\label{ss2.2}

The connection between the non-nullity of the powers of
Stiefel-Whitney characteristic classes and the lower bounds for graph
colorings is provided by the following general observation.

\begin{thm} \label{thmmain}
Let $G$ be a graph without loops, and let $T$ be a~graph with
$\zz$-action which flips some edge in $T$. If, for some integers
$k\geq 0$, $m\geq 1$, we have $\varpi_1^k(\thom(T,G))\neq 0$, and
$\varpi_1^{k}(\thom(T,K_m))=0$, then $\chi(G)\geq m+1$.
\end{thm}
\pr We have already shown that, under the assumptions of the theorem,
$\thom(T,H)$ is a~$\zz$-space for any loopfree graph $H$.  Assume now
that the graph $G$ is $m$-colorable, i.e., there exists a~homomorphism
$\phi:G\ra K_m$. It induces a~$\zz$-map
$\phi^T:\thom(T,G)\ra\thom(T,K_m)$. Since the Stiefel-Whitney classes
are functorial and $\varpi_1^{k}(\thom(T,K_m))=0$, the existence of
the $\zz$-map $\phi^{T}$ implies that $\varpi_1^{k}(\thom(T,G))=0$,
which is a~contradiction to the assumption of the theorem. \qed

\begin{rem} \label{rem:conn}
If a~$\zz$-space $X$ is $k$-connected, then there exists a~$\zz$-map
$\phi:S_a^{k+1}\ra X$, in particular, $\varpi_1^{k+1}(X)\neq 0$.
\end{rem}
\pr To construct $\phi$, subdivide $S_a^{k+1}$ simplicially as a join
of $k+2$ copies of $S^0$, and then define $\phi$ on the join of the
first $i$ factors, starting with $i=1$, and increasing $i$ by $1$ at
the time. To define $\phi$ on the first factor $\{a,b\}$, simply map
$a$ to an arbitrary point $x\in X$, and then map $b$ to $\gamma(x)$,
where $\gamma$ is the free involution of~$X$. Assume $\phi$ is defined
on $Y$ - the join of the first $i$ factors. Extend $\phi$ to
$Y*\{a,b\}$ by extending it first to $Y*\{a\}$, which we can do, since
$X$ is $k$-connected, and then extending $\phi$ to the second
hemisphere $Y*\{b\}$, by applying the involution~$\gamma$.

Since the Stiefel-Whitney classes are functorial, we have
$\phi^*(\varpi_1^{k+1}(X))=\varpi_1^{k+1}(S_a^{k+1})$, and the latter is
clearly nontrivial.\qed \vspace{5pt}

Let $T$ be any graph and consider the following equation
\begin{equation}\label{eq:ccvan}
\varpi_1^{n-\chi(T)+1}(\thom(T,K_n))=0, \text{ for all }
n\geq \chi(T)-1.
\end{equation}

\begin{thm} \label{thmmain1} $\,$
\begin{enumerate}
\item[(a)] The equation \eqref{eq:ccvan} is true for $T=K_m$, $m\geq 2$.
\item[(b)] The equation \eqref{eq:ccvan} is true for $T=C_{2r+1}$,
$r\geq 1$, and odd~$n$.
\end{enumerate}
\end{thm}
\pr The case $T=K_m$ is \cite[Theorem 1.6]{BK2} and has been proved
there. The case $T=C_{2r+1}$ will be proved in the
Section~\ref{sect5}. \qed

\begin{rem}
For a fixed value of $n$, if the equation \eqref{eq:ccvan} is true for
$T=C_{2r+1}$, then it is true for any $T=C_{2\ti r+1}$, if $r\geq\ti r$.
\end{rem}

\pr If $r\geq\ti r$, there exists a~graph homomorphism
$\phi:C_{2r+1}\ra C_{2\ti r+1}$ which respects the $\zz$-action.  This
induces a~$\zz$-map 
$$\phi_{K_n}:H^*(\thom(C_{2r+1},K_n))\ra H^*(\thom(C_{2\ti r+1},
K_n)),$$ 
yielding 
$$\ti\phi_{K_n}:H^*(\thom(C_{2r+1},K_n)/\zz;\zz)\ra H^*(\thom(C_{2\ti
r+1},K_n)/\zz;\zz).$$ Clearly,
$\ti\phi_{K_n}(\varpi_1(\thom(C_{2r+1},K_n)))= \varpi_1(\thom(C_{2\ti
r+1},K_n))$. In particular, $\varpi_1^i(\thom(C_{2r+1},K_n))=0$,
implies $\varpi_1^i(\thom(C_{2\ti r+1}, K_n))=0$.  \qed \vspace{5pt}

Note that for $T=C_{2r+1}$ and $n=2$, the equation \eqref{eq:ccvan} is
obvious, since $\thom(C_{2r+1},K_2)=\emptyset$. We give a~quick
argument for the next case $n=3$. One can see by inspection that the
connected components of $\thom(C_{2r+1},K_3)$ can be indexed by the
winding numbers $\alpha$. These numbers must be odd, so $\alpha=\pm
1,\pm 3,\dots,\pm(2s+1)$, where
\[
s=
\begin{cases}
(r-1)/3, & \text{ if } r\equiv 1 \mod 3,\\
\left\lfloor(r-2)/3\right\rfloor, & \text{ otherwise,}
\end{cases}
\] 
in particular $s\geq 0$. Let $\phi:\thom(C_{2r+1},K_3)\ra\{\pm 1,\pm
3, \dots,\pm(2s+1)\}$ map each point $x\in\thom(C_{2r+1},K_3)$ to the
point on the real line, indexing the connected component
of~$x$. Clearly, $\phi$ is a $\zz$-map.  Since $H^1(\{\pm 1,\pm 3,
\dots,\pm(2s+1)\}/\zz;\zz)=0$, the functoriality of the characteristic
classes implies $\varpi_1(\thom(C_{2r+1},K_3))=0$.

\begin{conj} \label{conj:sw}
The equation \eqref{eq:ccvan} is true for $T=C_{2r+1}$, $r\geq 1$, and
all~$n$.
\end{conj}


\subsection{Completing the sketch of the proof of the Lov\'asz Conjecture.} 
$\,$ \label{ss2.3} 

Consider one of the two maps $\iota:K_2\ra C_{2r+1}$ mapping the edge
to the $\zz$-invariant edge of $C_{2r+1}$. Clearly, $\iota$ is
$\zz$-equivariant. Since $\thom(-,H)$ is a~contravariant functor,
$\iota$ induces a~map of $\zz$-spaces
$\iota_{K_n}:\thom(C_{2r+1},K_n)\ra \thom(K_2,K_n)$, which in turn
induces a~$\dz$-algebra homomorphism
$\iota_{K_n}^*:H^*(\thom(K_2,K_n);\dz)\ra H^*(\thom(C_{2r+1},K_n);\dz)$.

\begin{thm}\label{thm:even_n}
Assume $n$ is even, then $2\cdot\iota_{K_n}^*$ is a $0$-map.
\end{thm}

Theorem \ref{thm:even_n} is proved in Section~\ref{sect4}. The results
of this paper were announced in \cite{BK1}, and the preprint of this
paper has been available since February 2004. In the summer 2005 an
alternative proof of Theorem~\ref{thm:even_n} appeared in the preprint
\cite{Z1}, and a~proof of Conjecture~\ref{conj:sw} was announced by
C.~Schultz.


\vspace{5pt}

\noindent
{\bf Proof of Theorem~\ref{loconj} (Lov\'asz Conjecture).}

\noindent
The case $k=-1$ is trivial, so take $k\geq 0$. Assume first that $k$
is even. By the Remark~\ref{rem:conn}, we have
$\varpi_1^{k+1}(\thom(C_{2r+1},G))\neq 0$. By
Theorem~\ref{thmmain1}(b), we have
$\varpi_1^{k+1}(\thom(C_{2r+1},K_{k+3}))=0$. Hence, applying
Theorem~\ref{thmmain} for $T=C_{2r+1}$ we get $\chi(G)\geq k+4$.

Assume now that $k$ is odd, and that $\chi(G)\leq k+3$. Let $\phi:G\ra
K_{k+3}$ be a~vertex-coloring map. Combining the
Remark~\ref{rem:conn}, the fact that $\thom(C_{2r+1},-)$ is
a~covariant functor from loopfree graphs to $\zz$-spaces, and the map
$\iota:K_2\ra C_{2r+1}$, we get the following diagram of $\zz$-spaces
and $\zz$-maps:
\begin{multline*}
S_a^{k+1}\stackrel{f}{\lra}\thom(C_{2r+1},G)
\stackrel{\phi^{C_{2r+1}}}{\lra}\thom(C_{2r+1},K_{k+3})
\stackrel{\iota_{K_{k+3}}}{\lra}\\
\stackrel{\iota_{K_{k+3}}}{\lra}\thom(K_2,K_{k+3})\cong S_a^{k+1}.
\end{multline*}
This gives a~homomorphism on the corresponding cohomology groups in
dimension $k+1$, $h^*=f^*\circ(\phi^{C_{2r+1}})^*
\circ(\iota_{K_{k+3}})^*:\dz\ra\dz$. It is well-known, see, e.g.,
\cite[Proposition 2B.6, p.\ 174]{Hat}, that a $\zz$-map $S_a^n\ra
S_a^n$ cannot induce a~$0$-map on the $n$th cohomology groups (in fact
it must be of odd degree). Hence, we have a~contradiction, and so
$\chi(G)\geq k+4$.  \qed \vspace{5pt}

Let us make a~couple of remarks.

\begin{rem}
As is apparent from our argument, we are actually proving a~sharper
statement than the original Lov\'asz Conjecture. First of all, the
condition ``$\thom(C_{2r+1},G)$ is $k$-connected'' can be replaced by
a~weaker condition ``the coindex of $\thom(C_{2r+1},G)$ is at least
$k+1$''. Furthermore, for even~$k$, that condition can be weakened
even further to ``$\varpi_1^{k+1}(\thom(C_{2r+1},G))\neq
0$''. Conjecture~\ref{conj:sw} would imply that this weakening can be
done for odd $k$ as well.
\end{rem}

\begin{rem}
It follows from \cite[Proposition 5.1]{BK2} that Lov\'asz Conjecture
is true if $C_{2r+1}$ is replaced by any graph $T$, such that $T$ can
be reduced to $C_{2r+1}$, by a~sequence of folds.
\end{rem}

\section{$\thomp$ and filtrations}

\subsection{The $+$ construction.} $\,$ \vspace{5pt}

For a finite graph $H$, let $H_+$ be the graph obtained from $H$ by
adding an extra vertex $b$, called the base vertex, and connecting it
by edges to all the vertices of $H_+$ including itself, i.e.,
$V(H_+)=V(H)\cup\{b\}$, and $E(H_+)=E(H)\cup\{(v,b),(b,v)\,|\,v\in
V(H_+)\}$.

\begin{df}
  Let $G$ and $H$ be two graphs. The simplicial complex $\thomp(G,H)$
  is defined to be the link in $\thom(G,H_+)$ of the homomorphism
  mapping every vertex of $G$ to the base vertex in $H_+$.
\end{df}

So the cells in $\thomp(G,H)$ are indexed by all $\eta:V(G)\ra
2^{V(H)}$ satisfying the same condition as in the
Definition~\ref{dfhom}. The closure of $\eta$ is also defined
identical to how it was defined for $\thom$. Note, that $\thomp(G,H)$
is simplicial, and that $\thomp(G,-)$ is a~covariant functor from {\bf
  Graphs} to {\bf Top}. One can think of $\thomp(G,H)$ as a~cell
structure imposed on the set of all partial homomorphisms from $G$
to~$H$.


\begin{figure}[hbt]
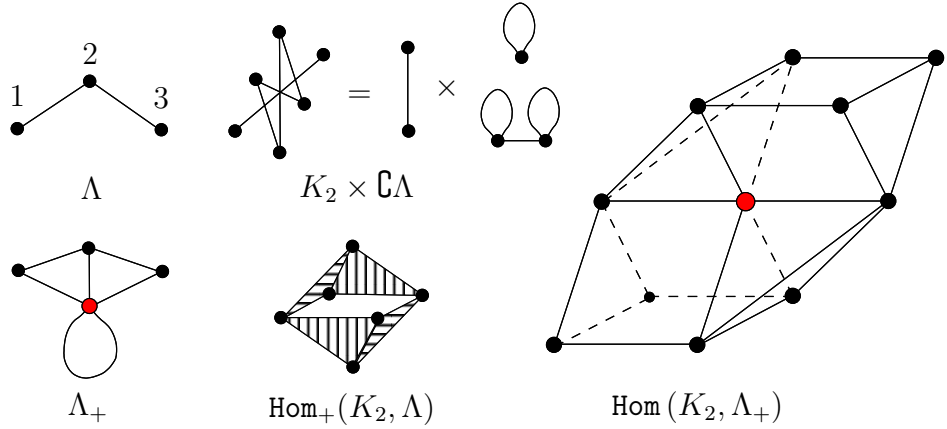

\begin{center}
  \begin{picture}(0,0)%
    \includegraphics{ex+.pstex}%
  \end{picture}%
  \input{ex+.pstex_t}%
 
\end{center}
\caption{The hom plus construction.}
\label{fig:ex+}
\end{figure}

For an arbitrary graph $G$, let $\ind(G)$ denote the independence
complex of $G$, i.e., the vertices of $\ind(G)$ are all vertices
of~$G$, and simplices are all the independent sets of~$G$.
The~dimension of $\thomp(G,H)$, unlike that of $\thom(G,H)$ is easy to
find:
\[
\dim(\thomp(G,H))=|V(H)|\cdot(\dim \ind(G)+1)-1.
\]

Recall that for any graph $G$, the strong complement $\com G$ is
defined by $V(\com G)=V(G)$, $E(\com G)=V(G)\times V(G)\sm
E(G)$. Also, for any two graphs $G$ and $H$, the direct product
$G\times H$ is defined by $V(G\times H)=V(G)\times V(H)$, $E(G\times
H)=\{((x,y),(x',y'))\,|\,(x,x')\in E(G),(y,y')\in E(H)\}$.

Sometimes, it is convenient to view $\thomp(G,H)$ as an~independence
complex of a~certain graph.

\begin{prop} \label{pr_homp}
The complex $\thomp(G,H)$ is isomorphic to $\ind(G\times\com H)$. In
particular, $\thomp(G,K_n)$ is isomorphic to $\ind(G)^{*n}$, where $*$
denotes the simplicial join.
\end{prop}

\pr By the definition, $V(G\times\com H)=V(G)\times V(H)$. Let
$S\subseteq V(G)\times V(H)$, $S=\{(x_i,y_i)\,|\,i\in I, x_i\in V(G),
y_i\in V(H)\}$. Then $S\in\ind(G\times\com H)$ if and only if, for any
$i,j\in I$, we have either $(x_i,x_j)\notin E(G)$ or $(y_i,y_j)\in
E(H)$, since the forbidden constellation is when $(x_i,x_j)\in E(G)$
and $(y_i,y_j)\notin E(H)$.

Identify $S$ with $\eta_S:V(G)\ra 2^{V(H)}$ defined by: for $v\in
V(G)$, set $\eta_S(v):=\{w\in V(H)\,|\,(v,w)\in S\}$. The condition
for $\eta_S\in\thomp(G,H)$ is that, if $(v_1,v_2)\in E(G)$, and
$w_1\in\eta_S(v_1)$, $w_2\in\eta_S(v_2)$, then $(w_1,w_2)\in E(H)$,
which is visibly identical to the condition for $S\in\ind(G\times\com
H)$. Hence $\thomp(G,H)=\ind(G\times\com H)$.

To see the second statement note first that $\com K_n$ is the
disjoint union of $n$ looped vertices. Since taking direct
products is distributive with respect to disjoint unions, and a
direct product of $G$ with a loop is again $G$, we see that
$G\times\com K_n$ is a disjoint union of $n$ copies of $G$.
Clearly, its independence complex is precisely the $n$-fold join
of $\ind(G)$. \qed

\subsection{Cochain complexes for $\thom(G,H)$ and $\thomp(G,H)$.} 
$\,$ \vspace{5pt} \label{ss3.2}

For any CW complex $K$, let $K^{(i)}$ denote the $i$-th skeleton
of~$K$. Let $R$ be a~commutative ring with a~unit. In this paper we
will have two cases: $R=\dz$ and $R=\zz$. For any $\eta\in K^{(i)}$,
we fix an~orientation on $\eta$, and let
$C_i(K;R):=R[\eta\,|\,\eta\in K^{(i)}]$, where
$R[\alpha\,|\,\alpha\in I]$ denotes the free $R$-module generated
by $\alpha\in I$. Furthermore, let $C^i(K;R)$ be the dual
$R$-module to $C_i(K;R)$. For arbitrary $\alpha\in C_i(K;R)$ let
$\alpha^*$ denote the element of $C^i(K;R)$ which is dual
to~$\alpha$. Clearly, $C^i(K;R)=R[\eta^*\,|\,\eta\in K^{(i)}]$,
and the cochain complex of $K$ is
\[\dots\stackrel{\bo^{i-1}}\lra C^i(K;R)\stackrel{\bo^i}\lra
C^{i+1}(K;R)\stackrel{\bo^{i+1}}\lra\dots.\] 

For $\eta\in K^{(i)}$, $\ti\eta\in K^{(i+1)}$, we have the incidence
number $[\eta:\ti\eta]$, which is $0$ if $\eta\notin\ti\eta$. In these
notations $\bo^i(\eta^*)=\sum_{\ti\eta\in K^{(i+1)}}$
$[\eta:\ti\eta]\,\ti\eta^*$. For arbitrary $\alpha\in C_i(K;R)$,
resp.\ $\alpha^*\in C^i(K;R)$, we let $[\alpha]$, resp.\ 
$[\alpha^*]$, denote the corresponding element of $H_i(K;R)$, resp.\ 
$H^i(K;R)$.

When coming after the name of a~cochain complex, the brackets $[-]$
will denote the index shifting (to the left), that is for the cochain
complex $C^*$, the cochain complex $C^*[s]$ is defined by
$C^i[s]:=C^{i+s}$, and the differential is the same (we choose not to
change the sign of the differential).

We now return to our context. Let $G$ and $H$ be two graphs, and let
us choose some orders on $V(G)=\{v_1,\dots,v_{|V(G)|}\}$ and on
$V(H)=\{w_1,\dots,w_{|V(H)|}\}$. Through the end of this subsection we
assume the coefficient ring to be $\dz$; the situation over $\zz$ is
simpler and can be described by tensoring with $\zz$.

Vertices of $\thomp(G,H)$ are indexed with pairs $(x,y)$, where $x\in
V(G)$, $y\in V(H)$, such that if $x$ is looped, then so is~$y$. We
order these pairs lexicographically:
$(v_{i_1},w_{j_1})\prec(v_{i_2},w_{j_2})$ if either $i_1<i_2$, or
$i_1=i_2$ and $j_1<j_2$. Orient each simplex of $\thomp(G,H)$
according to this order on the vertices. We call this orientation {\it
standard}, and call the oriented simplex~$\eta_+$. If
$\ti\eta_+\in\thomp^{(i+1)}(G,H)$ is obtained from
$\eta_+\in\thomp^{(i)}(G,H)$ by adding a~vertex~$v$, then
$[\eta_+:\ti\eta_+]$ is $(-1)^{k-1}$, where $k$ is the position of $v$
in the order of the vertices of~$\ti\eta_+$.

Let us now turn to $C^*(\thom(G,H))$. We can fix an orientation, which
we also call standard, on each cell $\eta\in\thom(G,H)$ as follows:
orient each simplex $\eta(i)$ according to the chosen order on the
vertices of $H$, then, order these simplices in the direct product
according to the chosen order on the vertices of~$G$. To simplify our
notations, we still call this oriented cell $\eta$, even though
a~choice of orders on the vertex sets of $G$ and $H$ is implicitly
present.

We remark for later use, that permuting the vertices of the simplex
$\eta(i)$ by some $\sigma\in\cs_{|\eta(i)|}$ changes the orientation
of the cell $\eta$ by $\sgn(\sigma)$, whereas swapping the simplices
with vertex sets $\eta(i)$ and $\eta(i+1)$ in the direct product
changes the orientation by
$(-1)^{(|\eta(i)|-1)(|\eta(i+1)|-1)}=(-1)^{\dim\eta(i)\cdot\dim\eta(i+1)}$.

If $\ti\eta\in\thom^{(i+1)}(G,H)$ is obtained from
$\eta\in\thom^{(i)}(G,H)$ by adding a~vertex~$v$ to the list
$\eta(t)$, then $[\eta:\ti\eta]$ is $(-1)^{k+d-1}$, where $k$ is the
position of $v$ in $\ti\eta(t)$, and $d$ is the dimension of the
product of the simplices with the vertex sets
$\eta(1),\dots,\eta(t-1)$, that is $d=1-t+\sum_{j=1}^{t-1}|\eta(j)|$.
To see this, note that $[\eta:\ti\eta]=1$ if the first vertex in the
first simplex is inserted. The general case follows from the
previously described rules for changing the sign of the orientation
under permuting simplices in the product and permuting vertices within
simplices.

\subsection{The support map and the relation between $\thom(G,H)$ 
and $\thomp(G,H)$.} \label{ss3.3a} $\,$ \vspace{5pt}

For each simplex of $\thomp(G,H)$, $\eta:V(G)\ra 2^{V(H)}$, define the
support of $\eta$ to be
$\supp\eta:=V(G)\setminus\eta^{-1}(\emptyset)$. A~concise way to
phrase the definition of $\supp$ differently is to consider the map
$t^G:\thomp(G,H)\ra\thomp(G,\com K_1)\simeq\Delta_{|V(G)|-1}$ induced
by the homomorphism $t:H\ra\com K_1$. Then, for each
$\eta\in\thomp(G,H)$ we have $\supp\eta=t^G(\eta)$, where the
simplices in $\Delta_{|V(G)|-1}$ are identified with the subsets of
$V(G)$.

Let $\wti C^*$ be the subcomplex of $C^*(\thomp(G,H))$ generated by
all $\eta^*_+$, for $\eta:V(G)\ra 2^{V(H)}$, such that
$\supp\eta=V(G)$ (cf.\ filtration in subsection \ref{ss3.3}). Set
\begin{equation} \label{eq:x*}
X^*(G,H):=\wti C^*[|V(G)|-1].
\end{equation}
Note that both $C^i(\thom(G,H))$ and $X^i(G,H)$ are free $\dz$-modules
with the bases $\{\eta^*\}_\eta$ and $\{\eta^*_+\}_\eta$ indexed by
$\eta:V(G)\ra 2^{V(H)}\sm\{\emptyset\}$, such that
$\sum_{j=1}^{|V(G)|}|\eta(j)|=|V(G)|+i$.

At this point we introduce the following notations: for $\eta:V(G)\ra
2^{V(H)}\sm\{\emptyset\}$, set 
$$c(\eta):=\sum_{\begin{subarray}{c}i\text{ is even}\\
1\leq i\leq|V(G)|\end{subarray}}
|\eta(i)|.$$

For any $\eta:V(G)\ra 2^{V(H)}\sm\{\emptyset\}$, set
$\rho(\eta_+):=(-1)^{c(\eta)}\eta$. Obviously, the induced map
$\rho^*:X^i(G,H)\ra C^i(\thom(G,H))$ is a~$\dz$-module isomorphism for
any~$i$.

\begin{prop}
 The map $\rho^*:X^*(G,H)\ra C^*(\thom(G,H))$ is an~isomorphism of the
  cochain complexes.
\end{prop}

\pr Indeed, let $\ti\eta:V(G)\ra 2^{V(H)}\sm\{\emptyset\}$ be obtained
from $\eta$ by adding a~vertex $v$ to the list $\eta(t)$, and let $k$
be the position of $v$ in $\ti\eta(t)$. By our previous computation:
$[\eta_{+}:\ti\eta_+]=(-1)^{k+d+t}$, whereas
$[\eta:\ti\eta]=(-1)^{k+d-1}$, where
$d=1-t+\sum_{j=1}^{t-1}|\eta(j)|$. This shows that
$$
[\rho(\eta_+):\rho(\ti\eta_+)]=(-1)^{c(\eta)+c(\ti\eta)}
[\eta:\ti\eta]= (-1)^{c(\eta)+c(\ti\eta)+t+1}[\eta_+:\ti\eta_+],
$$
but
$$c(\eta)+c(\ti\eta)+t+1=\sum_{i\text{ is even}}|\eta(i)|+
\sum_{i\text{ is even}}|\ti\eta(i)|+t+1\equiv 0\,\,(\text{mod } 2)$$
for~any~$t$, hence $[\rho(\eta_+):\rho(\ti\eta_+)]=
[\eta_+:\ti\eta_+]$. \qed

\subsection{Relating $\zz$-actions on $\thom(G,H)$ and $\thomp(G,H)$.}  
$\,$ \vspace{5pt}

Assume that we have $\gamma\in\chom(G,G)$, and $0\leq r\leq |V(G)|/2$,
such that
\[\gamma(i)=\begin{cases}
2r+1-i,& \text{ if }1\leq i\leq 2r; \\
i,&\text{ if }2r+1\leq i\leq |V(G)|,
\end{cases}\]
where we identified $V(G)$ with the numbering $[1,V(G)]$. In
particular, we have $\gamma^2=1$.

The homomorphism $\gamma$ induces $\zz$-action both on $\thom(G,H)$,
and on $\thomp(G,H)$. We shall see how $\rho^*$ behaves with respect
to this $\zz$-action. For any $\eta:V(G)\ra 2^{V(H)}\sm\{\emptyset\}$,
$\gamma$ takes $\eta$ to $\eta\circ\gamma$. By a~slight abuse of
notations we let $\gamma$ denote the induced actions both on
$C^*(\thom(G,H))$ and on $X^*(G,H)$.

Let $(u_1,\dots,u_q)$ be the vertices of the simplex $\eta_+$ listed
in the increasing order. By definition,
$\gamma(\eta_+)=(\gamma(u_1),\dots,\gamma(u_q))$, where
$\gamma((v,w)):=(\gamma(v),w)$, for $v\in V(G)$, $w\in V(H)$. Clearly,
$\gamma(\eta_+)$ has the same set of vertices as
$(\eta\circ\gamma)_+$, so we just need to see how their orientations
relate. To order the vertices of $\gamma(\eta_+)$ we need to invert
the order of the blocks with cardinalities
$|\eta(1)|,\dots,|\eta(2r)|$ without changing the vertex orders within
the blocks. The sign of this permutation is $(-1)^c$, where
$c=\sum_{1\leq i<j\leq 2r}|\eta(i)|\cdot|\eta(j)|$, so we conclude
that 
\begin{equation}
  \label{eq:geta+}\gamma(\eta^*_+)=(-1)^c(\eta\circ\gamma)^*_+.
\end{equation}

Consider now the oriented cell $\eta$. It is a~direct product of
simplices $\da^1,\dots,\da^{|V(G)|}$ of dimensions
$|\eta(1)|-1,\dots,|\eta(|V(G)|)|-1$, with the standard orientation as
defined above. The cell $\gamma(\eta)$ is the direct product of
$\gamma(\da^1)=\da^{2r},\gamma(\da^2)=\da^{2r-1},\dots,
\gamma(\da^{2r})=\da^1,\gamma(\da^{2r+1})=\da^{2r+1},$ $\dots,
\gamma(\da^{|V(G)|})=\da^{|V(G)|}$, with the order of the vertices
(hence the orientation) within each simplex being the same as
in~$\eta$.

We see that $\gamma(\eta)$ is, up to the orientation, the same cell as
$\eta\circ\gamma$. To relate their orientations, we need to permute
the simplices $\da^{2r},\dots,\da^1$ back in order, which, by the
previous observations, changes the orientation by $(-1)^{\ti d}$,
where
\begin{multline*}
\ti d=\sum_{1\leq i<j\leq 2r}\dim\da^i\cdot\dim\da^j=
\sum_{1\leq i<j\leq 2r}(|\eta(i)|-1)(|\eta(j)|-1)\\
=c-(2r-1)\sum_{i=1}^{2r}|\eta(i)|+{\binom{2r}{2}}.
\end{multline*}
Reducing modulo 2, we conclude that 
\begin{equation}
  \label{eq:geta}
  \gamma(\eta^*)=(-1)^d(\eta\circ\gamma)^*,
\end{equation}
where $d=c+\sum_{i=1}^{2r}|\eta(i)|+r$.

Let us now see how $\rho^*$ interacts with $\gamma$. We have
\begin{equation}
  \label{eq:3.3}
  \rho^*(\gamma(\eta^*_+))=(-1)^c\rho^*((\eta\circ\gamma)^*_+)=
(-1)^{c+c(\eta\circ\gamma)}(\eta\circ\gamma)^*,
\end{equation}
where the first equality is by \eqref{eq:geta+} and the second one is
by definition of $\rho$, and
\begin{equation}
  \label{eq:3.4}
  \gamma(\rho^*(\eta^*_+))=(-1)^{c(\eta)}\gamma(\eta^*)=
(-1)^{d+c(\eta)}(\eta\circ\gamma)^*,
\end{equation}
where the first equality is by definition of $\rho$ and second one is
by~\eqref{eq:geta}. Comparing \eqref{eq:3.3} with \eqref{eq:3.4}, and
using the computation
\begin{multline*}
c(\eta)+c(\eta\circ\gamma)=\sum_{\begin{subarray}{c}{i\text{ is even}}\\
{1\leq i\leq|V(G)|}\end{subarray}}|\eta(i)|+
\sum_{\begin{subarray}{c}{i\text{ is even}}\\{1\leq i\leq|V(G)|}
\end{subarray}}|\eta\circ\gamma(i)|\\
=\sum_{i=1}^{2r}|\eta(i)|+2\cdot\sum_{\begin{subarray}{c}{i\text{ is even}}\\
{2r+1\leq i\leq|V(G)|}\end{subarray}}|\eta(i)|,
\end{multline*}
we see that, for any $\eta$
\begin{equation}
  \label{eq:act_disc}
   \rho^*(\gamma(\eta^*_+))=(-1)^r \gamma(\rho^*(\eta^*_+)).
\end{equation}

\subsection{The filtration of $C^*(\thomp(G,H);\dz)$ and the 
  $E_0^{*,*}$-tableau.} \label{ss3.3} $\,$ \vspace{5pt}

We shall now filter $C^*(\thomp(G,H);\dz)$. Define the subcomplexes of
$C^*(\thomp(G,H);\dz)$, $F^p=F^pC^*(\thomp(G,H);\dz)$, as follows:
\[ F^p: \dots\stackrel{\bo^{q-1}}\lra F^{p,q}\stackrel{\bo^{q}}\lra
F^{p,q+1}\stackrel{\bo^{q+1}}\lra\dots,\]
where
\[ F^{p,q}=F^p C^q(\thomp(G,H);\dz)=\dz\left[\eta^*_+\,\left|\,
\eta_+\in\thomp^{(q)}(G,H),|\supp\eta|\geq p+1\right]\right.,
\]
and $\bo^*$ is the restriction of the differential in
$C^*(\thomp(G,H);\dz)$. Then,
\[
C^q(\thomp(G,H);\dz)=F^{0,q}\supseteq F^{1,q}\supseteq\dots\supseteq
F^{|V(G)|-1,q}\supseteq F^{|V(G)|,q}=0.
\]

\begin{prop}
For any $p$,
\begin{equation}\label{eq:e0tab}
F^p/F^{p+1}=\bigoplus_{\begin{subarray}{c}{S\subseteq V(G)}\\{|S|=p+1}
\end{subarray}}
C^*(\thom(G[S],H);\dz)[-p].
\end{equation}
Hence, the 0th tableau of the spectral sequence associated to
the cochain complex filtration $F^*$ is given by
\begin{equation}
  \label{eq:E0}
E_0^{p,q}=C^{p+q}(F^p,F^{p+1})=\bigoplus_{\begin{subarray}{c}
{S\subseteq V(G)}\\{|S|=p+1}\end{subarray}}
C^q(\thom(G[S],H);\dz).
\end{equation}
\end{prop}
\pr By construction
\begin{multline*}
  F^{p,q}/F^{p+1,q}=\dz\left[\eta^*_+\,\left|\,\eta_+\in\thomp^{(q)}(G,H),
  |\supp\eta|=p+1\right]\right.=\\
  \bigoplus_{\begin{subarray}{c}{S\subseteq V(G)}\\{|S|=p+1}\end{subarray}}
X^{q-p}(G[S],H;\dz)
  \stackrel{\rho^*}{=} 
\bigoplus_{\begin{subarray}{c}{S\subseteq V(G)}\\{|S|=p+1}\end{subarray}}
C^{q-p}(\thom(G[S],H);\dz),
\end{multline*} 
where $X^*$ is defined in \eqref{eq:x*}, and $\rho^*$ is the map
defined in subsection~\ref{ss3.3a}.  \qed \vspace{5pt}

Note, that in particular we have
$F^{|V(G)|-1,q}/F^{|V(G)|,q}=F^{|V(G)|-1,q}=C^q(\thom(G,H);\dz)[1-|V(G)|]$.

\section{$\zz$-action on $H^*(\thomp(C_{2r+1},K_n);\dz)$}
\label{sect4}

In this section we shall derive some information about the
$\dz[\zz]$-modules $H^*(\thom(C_{2r+1},K_n);\dz)$, for $r\geq 2$,
$n\geq 4$. For $r=2$ our computation will be complete.

We adopt the following convention: we think of $C_{2r+1}$ as a~unit
circle on the plane with $2r+1$ marked points with numbers
$1,\dots,2r+1$ following each other in the clockwise increasing order.
The directions {\it left}, resp.\ {\it right} on this circle will
denote counterclockwise, resp.\ clockwise. 

Furthermore, before we start our computation, we introduce the
following terminology. For $S\subset V(C_{2r+1})$, we call those
connected components of $C_{2r+1}[S]$ which have at least 2 vertices
the {\it arcs}.  For $x,y\in\dz$, we let $[x,y]_{2r+1}$ denote the arc
starting from $x$ and going clockwise to $y$, that is
$[x,y]_{2r+1}=\{[x]_{2r+1},[x+1]_{2r+1},\dots,[y-1]_{2r+1},[y]_{2r+1}\}$.

\subsection{The simplicial complex of partial homomorphisms from a~cycle
  to a~complete graph.} $\,$ \vspace{5pt}

Here and in the next subsection we summarize some previously published
results which are necessary for our present computations. To start
with, recall that the homotopy type of the independence complexes of
cycles was computed in~\cite{Koz1}.

\begin{prop}\label{pr_indcomp}{\rm (\cite[Proposition 5.2]{Koz1})}. 

\nin For any integer $m\geq 2$, we have
\[\ind(C_m)\simeq\begin{cases}
S^{k-1}\vee S^{k-1},&\text{ if } m=3k;\\
S^{k-1},&\text{ if } m=3k\pm 1.
\end{cases}\]
\end{prop}

Combining Propositions~\ref{pr_homp} and \ref{pr_indcomp} we get the
following formula.

\begin{crl}\label{crl_indcomp} 
 For any integers $m\geq 2$, $n\geq 1$, we have
\begin{equation}\label{eq_indcomp}
\thomp(C_m,K_n)\simeq\begin{cases}
\bigvee_{2^n\text{ copies}} S^{nk-1},&\text{ if } m=3k;\\
S^{nk-1},&\text{ if } m=3k\pm 1.
\end{cases}
\end{equation}
\end{crl}
 
The following estimates will be needed later for our spectral sequence
computations.

\begin{crl}\label{crl:h+}
We have $\wti H^i(\thomp(C_{2r+1},K_n))=0$ for $r\geq 2$, $n\geq 4$,
and $i\leq n+2r-2$; except for the two cases $(n,r)=(4,3)$ and
$(5,3)$.
\end{crl}
\pr Note, that if $2r+1=3k+\epsilon$, with $\epsilon\in\{-1,0,1\}$,
then we have $\wti H^i(\thomp(C_{2r+1},K_n))=0$, for $i\leq nk-2$.

Assume first $2r+1=3k$. The inequality $nk-2\geq n+2r-2$ is equivalent
to $n\geq 3+2/(k-1)$, and the latter is always true since $k\geq 3$
and $n\geq 4$.

 Assume now $2r+1=3k+1$. This time, $nk-2\geq n+2r-2$ is equivalent to
$n\geq 3+3/(k-1)$. If $k\geq 4$, this is always true, since $n\geq 4$.
If $k=2$, this reduces to saying that $n\geq 6$. This yields the two
exceptional cases: $r=3$ and $n=4,5$.

Finally, assume $2r+1=3k-1$. Here, $nk-2\geq n+2r-2$ is equivalent to
$n\geq 3+1/(k-1)$, which is always true, since $k\geq 2$, $n\geq 4$.
\qed \vspace{5pt}

The Corollary \ref{crl_indcomp} can be strengthened to include the
information on the $\zz$-action.

\begin{prop}\label{pr:h+/z2}
For any positive integers $r$ and $n$ we have
\begin{equation}\label{eq:h+/z2}
\thomp(C_{2r+1},K_n)/\zz\simeq
\begin{cases}
\bigvee_{2^{n-1}\text{ copies}} S^{nk-1},&\text{ if } 2r+1=3k;\\
S^{kn/2-1}*\rp^{kn/2-1},&\text{ if } 2r+1=3k\pm 1.
\end{cases}
\end{equation}
\end{prop}
\pr
By Proposition~\ref{pr_homp} we know that $\thomp(C_{2r+1},K_n)$
is isomorphic to $\ind(C_{2r+1})^{*n}$. We analyze $\zz$-action on
$\ind(C_{2r+1})$ in more detail.

Assume first $2r+1=3k-1$, in particular, $k$ is even. It was shown
in~\cite[Proposition~5.2]{Koz1} that
$X=\ind(C_{2r+1})\sm\{1,4,\dots,2r-3,2r\}$ is contractible (here
``$\sm$'' just means the removal of an open maximal simplex).  It
follows from the standard fact in the theory of transformation groups,
see e.g., \cite[Theorem 5.16, p.\ 222]{TtD}, that $X/\zz$ is
contractible as well. Hence $\ind(C_{2r+1})$ is $\zz$-homotopy
equivalent to the unit sphere $S^{k-1}\subset \dr^k$ with the $\zz$
acting by fixing $k/2$ coordinates and multiplying the other $k/2$
coordinates by $-1$.

Assume $2r+1=3k+1$. The link of the vertex $2r+1$ is $\zz$-homotopy
equivalent to a~point. Hence, deleting the open star of the vertex
$2r+1$ produces a~complex $X$, which is $\zz$-homotopy equivalent to
$\ind(C_{2r+1})$.  It was shown in \cite[Proposition~5.2]{Koz1}, that
$X\sm\{2,5,\dots,2r-4,2r-1\}$ is contractible.  By an~argument,
similar to the previous case, we conclude that $\ind(C_{2r+1})$ is
$\zz$-homotopy equivalent to the unit sphere $S^{k-1}\subset \dr^k$
with the $\zz$ acting by fixing $k/2$ coordinates and multiplying the
other $k/2$ coordinates by $-1$.

In both cases we see that $\thomp(C_{2r+1},K_n)$ is $\zz$-homotopy
equivalent to $\susp^{kn/2}S^{kn/2-1}$, with the $\zz$-action and the
latter space being induced by the antipodal action on $S^{kn/2-1}$. It
follows that $\thomp(C_{2r+1},K_n)/\zz$ is homotopy equivalent to
$\susp^{kn/2}\rp^{kn/2-1}$.

Consider the remaining case $2r+1=3k$. It was shown in
\cite[Proposition~5.2]{Koz1} that $\ind(C_{2r+1})$ becomes
contractible if one removes the simplices $\{1,4,\dots,2r-1\}$ and
$\{2,5,\dots,2r\}$. It follows that $\ind(C_{2r+1})$ is $\zz$-homotopy
equivalent to the wedge of two unit spheres $S^{k-1}$ with the $\zz$
acting by swapping the spheres. Thus $\thomp(C_{2r+1},K_n)$ is
$\zz$-homotopy equivalent to a~wedge of $2^n$ $(nk-1)$-dimensional
spheres, with the $\zz$-action swapping them in pairs, and so
$\thomp(C_{2r+1},K_n)/\zz$ is homotopy equivalent to a~wedge of
$2^{n-1}$ $(nk-1)$-dimensional spheres.
\qed \vspace{5pt}

We summarize the estimates which we will need later.

\begin{crl}
We have $\wti H^i(\thomp(C_{2r+1},K_n)/\zz)=0$ for $r\geq 2$, $n\geq
5$, and $i\leq n+r-2$. Except for the case $r=3$.
\end{crl}

\pr If $2r+1=3k$, the inequality $nk-2\geq n+r-2$ is equivalent to
$n\geq 3r/(2r-2)$, which is true for $n\geq 3$, $r\geq 2$. If
$2r+1=3k-1$, then $nk/2\geq n+r-2$ is equivalent to $(n-3)(k-2)\geq
0$, again true for our parameters.

If $2r+1=3k+1$, then $nk/2\geq n+r-2$ is equivalent to $(n-3)(k-2)\geq
2$. This is true for all parameters $n\geq 5$, $k\geq 2$, except for
$k=2$.  
\qed

\subsection{The cell complex of homomorphisms from a~tree
  to a~complete graph.} $\,$ 
\label{ss4.2}

In the next proposition we summarize several results proved
in~\cite{BK2,K4}.
\begin{prop}\label{pr:BK2}\cite[Propositions 4.3, 5.4, and 5.5]{BK2}, 
\cite{K4}.

\nin Let $T$ be a~tree with at least one edge.

\nin {\rm (i)} The map $i_{K_n}:\thom(T,K_n)\ra\thom(K_2,K_n)$ induced
by any inclusion $i:K_2\hookrightarrow T$ is a~homotopy equivalence.

\nin {\rm (ii)} $\thom(K_2,K_n)$ as a boundary complex of a polytope
of dimension $n-2$, in particular $\thom(T,K_n)\simeq S^{n-2}$.

\nin {\rm (iii)} Given a~$\zz$-action determined by an~invertible
graph homomorphism $\gamma:T\ra T$, if $\gamma$ flips an~edge in $T$,
then $\thom(T,K_n)\simeq_\zz S^{n-2}_a$, otherwise
$\thom(T,K_n)\simeq_\zz S^{n-2}_t$.
\end{prop}

\nin Here $S^{m}_a$ denotes the $m$-sphere equipped with an antipodal
$\zz$-action, whereas $S^{m}_t$ is the $m$-sphere equipped with the
trivial one.

Let $F$ be any graph, with $F_1,\dots,F_t$ being the list of all those
connected components of $F$ which have at least 2 vertices. For any
$\emptyset\neq S\subseteq[1,t]$, and $V=\{v_i\}_{i\in S}$, such that
$v_i\in V(F_i)$, for any $i\in S$, set
$$\alpha_+(F,V):=\sum_{\eta}\eta^*_+,\quad\alpha(F,V):=\sum_{\eta}\eta^*,$$
where both sums are taken over all $\eta:V(F)\ra
2^{[1,n]}\sm\{\emptyset\}$, such that
\begin{itemize}
\item $\eta(v_i)=[1,n-1]$, for all $i\in S$;
\item $|\eta(w)|=1$, for all $w\in V(F)\sm V$.
\end{itemize}
Note that, for fixed $S$ and $V$, $(-1)^{c(\eta)}$ does not depend on
the choice of $\eta$ as long as $\eta$ satisfies these two
conditions. In our previous notations we have $\alpha_+(F,V)\in
X^{|S|(n-2)}(F,K_n)$, and $\alpha(F,V)\in
C^{|S|(n-2)}(\thom(F,K_n))$. When $|S|=1$, $V=\{v\}$, we shall simply
write $\alpha_+(F,v)$ and $\alpha(F,v)$.

Assume now that $F$ is a forest. For $w\in V(F_i)$, such that
$(v_i,w)\in E(F)$, set
$W:=\{v_1,\dots,v_{i-1},w,v_{i+1},\dots,v_t\}=(V\cup\{w\})\sm\{v_i\}$.
We have a~graph homomorphism $K_2\ra(v,w)$, which induces
a~$\zz$-equivariant map $\varphi^*:H^*(\thom(K_2,K_n))\ra
H^*(\thom(F,K_n))$.  We know that $\thom(K_2,K_n)\cong_\zz S^{n-2}_a$,
and that the dual of any $(n-2)$-dimensional cell of $\thom(K_2,K_n)$
is a~generator of $H^{n-2}(\thom(K_2,K_n);\dz)$. Comparing
orientations of the cells of $\thom(K_2,K_n)$ we see that
$[\alpha(K_2,1)]=(-1)^{n-1}[\alpha(K_2,2)]$, where $1$ and $2$ denote
the vertices of $K_2$.  Applying $\varphi^*$ we conclude that
\[ [\alpha(F,V)]=(-1)^{n-1}[\alpha(F,W)].\]
Since $\rho^*$
is a~cochain isomorphism and
$\rho^*(\alpha_+(F,V))=(-1)^{c(\eta)}\alpha(F,V)$, we have
\begin{equation} \label{eq:my7}
[\alpha_+(F,V)]=
  \begin{cases}
    -[\alpha_+(F,W)],& \text{if $v$ and $w$ have different} \\ 
                     & \text{parity in the order on $V(F)$}; \\
    (-1)^{n-1}[\alpha_+(F,W)],& \text{if they have the same parity.}
  \end{cases}
\end{equation}

\subsection{The $E_1^{*,*}$-tableau for $E_1^{p,q}\Rightarrow 
H^{p+q}(\thomp(C_{2r+1},K_n);\dz)$.}\footnote{The calculations performed in the subsections 
\ref{ssect:4.3}--\ref{ssect:4.8} have been verified and generalized 
in~\cite{ccc}.}  $\,$ \vspace{5pt} \label{ssect:4.3}

Let us fix integers $r\geq 2$ and $n\geq 4$. Let
$(F^p)_{p=0,\dots,|V(G)|-1}$ be the filtration on
$C^*(\thomp(C_{2r+1},K_n);\dz)$ defined in subsection~\ref{ss3.3}, and
consider the corresponding spectral sequence. The entries of the
$E_1$-tableau are given by $E_1^{p,q}=H^{p+q}(F^p,F^{p+1})$.  Since
all proper subgraphs of $C_{2r+1}$ are forests, we can now use the
formula~\eqref{eq:E0} to obtain almost complete information about the
$E_1$-tableau. See Figure~\ref{fig:1e1}, where all the entries outside
of the shaded area are equal to~0

Let $\emptyset\neq S\subset V(C_{2r+1})$, and let
$S_1,\dots,S_{l(S)}$ be the connected components of $C_{2r+1}[S]$,
with $|S_{1}|\geq|S_{2}|\geq\dots\geq |S_{d(S)}|>|S_{d(S)+1}|=
\dots=|S_{l(S)}|=1$, where $l(S)\geq 1$, but we may have $d(S)=0$ or
$d(S)=l(S)$. By Proposition~\ref{pr:BK2} together with property (3)
from \cite[subsection 2.4]{BK2} we see that
\begin{equation}\label{eq:4.3a}
\thom(C_{2r+1}[S],K_n)\simeq\prod_{i=1}^{d(S)}S^{n-2}.
\end{equation}
Combining with the formula~\eqref{eq:E0} we conclude
\begin{equation}
  \label{eq:4.3b}
  H^{p+q}(F^p,F^{p+1})=\bigoplus_{\begin{subarray}{c}
{S\subset V(C_{2r+1})}\\{|S|=p+1}\end{subarray}}
H^q\left(\prod_{i=1}^{d(S)}S^{n-2};\dz\right),
\end{equation}
for $p\leq 2r-1$.


\begin{figure}[hbt]
\begin{center}
  \begin{picture}(0,0)%
    \includegraphics{1e1.pstex}%
  \end{picture}%
  \input{1e1.pstex_t}%
 
\end{center}
\caption{The $E_1^{*,*}$-tableau, for 
$E_1^{p,q}\Rightarrow H^{p+q}(\thomp(C_{2r+1},K_n);\dz)$.}
\label{fig:1e1}
\end{figure}

Since the spectral sequence converges to $\wti H^*
(\thomp(C_{2r+1},K_n) ;\dz)$, and, since by Corollary~\ref{crl:h+},
$\wti H^{i} (\thomp(C_{2r+1},K_n);\dz)=0$ for $i\leq n+2r-2$, we know
that the entries on the diagonals $p+q=n+2r-2$, and $p+q=n+2r-3$,
should eventually all become~$0$.

\subsection{The cochain complex $(D_0^*,d_1)=(E_1^{*,0},d_1)$, for
$E_1^{p,q}\Rightarrow H^{p+q}(\thomp(C_{2r+1},K_n);\dz)$.}
$\,$ \vspace{5pt}

Let $(D_i^*,d_1)$ denote the cochain complex in the $i(n-2)$-th row of
$E_1^{*,*}$, for any
$i=0,\dots,\left\lfloor\dfrac{2r+1}{3}\right\rfloor$. Next we show
that $(D_0^*,d_1)$ is isomorphic to the cochain complex of a~simplex.

\begin{lm}
We have $E_2^{0,0}=\dz$, and
$E_2^{1,0}=E_2^{2,0}=\dots=E_2^{2r,0}=0$.
\end{lm}

\pr Let $\Delta_{2r}$ denote and abstract simplex with $2r+1$ vertices
indexed by $[1,2r+1]$, and identify simplices of $\Delta_{2r}$ with
the subsets of $[1,2r+1]$. Let $(C^*(\Delta_{2r};\dz),d^*)$ be the
cochain complex of $\Delta_{2r}$ corresponding to the order on the
vertices given by this indexing. By (\ref{eq:4.3b}), each $S\subseteq
V(C_{2r+1})$, $|S|=p+1$, contributes one independent generator
(over~$\dz$) to $E_1^{p,0}$. Identifying it with the generator in
$C^*(\Delta_{2r};\dz)$ of the corresponding $p$-simplex in $\da_{2r}$,
we see that $(D_0^*,d_1)$ and $(C^*(\Delta_{2r};\dz),d^*)$ are
isomorphic as cochain complexes.

Indeed, for such $S$,
$\tau_S:=\sum_{\varphi\in\chom(C_{2r+1}[S],K_n)}\varphi^*_+$ is
a~representative of the corresponding generator in $E_1^{p,0}$.  This
is true even for $S=V(C_{2r+1})$, since $\thom(C_{2r+1},K_n)$ is
connected for $n\geq 4$, as was shown in~\cite[Proposition 2.1]{BK2}.
Clearly, 
\begin{multline*}  
 d_1(\tau_S)=
\sum_{\varphi\in\chom(C_{2r+1}[S],K_n)}\sum_{v\notin S}
\sum_{\psi|S=\varphi}[\varphi_+:\psi_+]\psi^*_+=\\
\sum_{v\notin S}[S:S\cup\{v\}]
\sum_{\psi\in\chom(C_{2r+1}[S\cup\{v\}],K_n)}\psi^*_+=
\sum_{v\notin S}[S:S\cup\{v\}]\tau_{S\cup\{v\}},
\end{multline*}
where the second equality is true since $[\varphi_+:\psi_+]$ only
depends on $S$ and $v$, not on the specific choice of $\varphi$
and~$\psi$.  This shows that the following diagram commutes
\[\begin{CD}
C^p(\Delta_{2r};\dz)@>d^p>>C^{p+1}(\Delta_{2r};\dz) \\
        @V\tau_\cdot VV           @VV\tau_\cdot V\\
E_1^{p,0} @>d_1>>E_1^{p+1,0}, \\
\end{CD}\]
where $\tau_\cdot:C^*(\Delta_{2r};\dz)\ra E_1^{*,0}$ is the linear
extension of the map taking $S$ to $\tau_S$, for $S\subseteq
V(C_{2r+1})$.  It follows that $(D_0^*,d_1)$ is isomorphic to
$(C^*(\Delta_{2r};\dz),d^*)$, therefore $E_2^{0,0}=\dz$, and
$E_2^{1,0}=E_2^{2,0}=\dots=E_2^{2r,0}=0$.
\qed

\subsection{The cochain complexes $(D_t^*,d_1)=(E_1^{*,(n-2)t},d_1)$, 
for $\left\lfloor(2r+1)/3\right\rfloor\geq t\geq 2$, and
$E_1^{p,q}\Rightarrow H^{p+q}(\thomp(C_{2r+1},K_n);\dz)$.} 
$\,$ \vspace{5pt}
\label{ss4.5}

We shall perform only a~partial analysis of the cohomology groups
of $(D_t^*,d_1)$, which will however be sufficient for our purpose.

For $S\subset V(C_{2r+1})$ and $v\in S$, let $a(S,v)$ denote the arc
of $S$ to which $v$ belongs (assuming this arc exists). Furthermore,
for an~arbitrary arc $a$ of $S$, let $a=[a_\bu,a^\bu]_{2r+1}$. Let
$|a|$ denote the number of vertices on $a$, and set $\widehat a:=
[a_\bu-1,a^\bu+1]_{2r+1}$ (so $|\widehat a|=|a|+2$, if $|a|\leq 2r-1$).

For any $V\subseteq S\subseteq V(C_{2r+1})$, as in
subsection~\ref{ss4.2}, set
$\sigma_{S,V}:=\alpha_+(C_{2r+1}[S],V)$. By our previous observations,
$E_1^{0,t(n-2)}=0$. Furthermore, for any $1\leq i\leq 2r-1$,
$E_1^{i,t(n-2)}$ is a~free $\dz$-module with the basis
$\{[\sigma_{S,V}]\}$, where $S\subset V(C_{2r+1})$, $|S|=i+1$,
$|V|=t$, and $v=a_\bu(S,v)$ (i.e., $[v-1]_{2r+1}\notin S$) for all
$v\in V$. Since $\sigma_{S,v}$ is a~cocycle in
$X^{n-2}(C_{2r+1}[S],K_n;\dz)$, we have
\begin{equation} \label{eq:my6}
  d_1([\sigma_{S,v}])=\sum_{w\notin S}(-1)^{z(w)}[\sigma_{S\cup\{w\},v}],
\end{equation}
where
$$z(w)=
\begin{cases}
  \left|S\cap[1,w-1]\right|,& \text{if } v\notin[1,w-1];\\
n+\left|S\cap[1,w-1]\right|,& \text{if } v\in[1,w-1].
\end{cases}
$$ Note, that if $i\leq 2r-2$ and $[w]_{2r+1}=[v-1]_{2r+1}$, then
$v\neq a_\bu(S\cup\{w\},w)$, so $[\sigma_{S\cup\{w\},v}]$ may differ
by a~sign from one of the elements in our chosen basis.  We shall not
need the analog of the equation~\eqref{eq:my6} for the case $|V|\geq
2$.

Let $A^*_1$ be the subcomplex of $D_1^*$ defined by:
$$A^*_1:0\lra\wti E_1^{2r-2,n-2}{\stackrel{d_1}\lra} E_1^{2r-1,n-2}
{\stackrel{d_1}\lra} E_1^{2r,n-2}\lra 0,$$
where the $\dz$-modules
indexed with $0,\dots,2r-3$ are equal to~0, and $\wti E_1^{2r-2,n-2}$
is generated by $\{[\sigma_{S,v}]\}$, such that $S$ and $v$ satisfy
all the previously required conditions and, in addition, $C_{2r+1}[S]$
is connected. 

In general, let $A^*_t$ be the subcomplex of $D^*_t$ generated by all
$\{[\sigma_{S,V}]\}$, such that
\begin{equation}\label{eq:hat}
\bigcup_{v\in V}\widehat{a(S,v)}=V(C_{2r+1}).
\end{equation}
In words: the gaps between those arcs of $S$ which have points in $V$
are of length at most~$2$. For future reference, we note, that
\eqref{eq:hat} implies that $|S|+2|V|\geq 2r+1$, i.e., $|S|-1\geq
2r-2t$, hence $A^j_t=0$ for $j<2r-2t$.

\begin{lm}
We have $H^*(D_t^*)=H^*(A_t^*)$.
\end{lm}

\pr
Let us set up another spectral sequence for computing the cohomology
of the relative complex $(D_t^*,A^*_t)$. We filter by $\sum_{v\in V}
|a(S,v)|$. More precisely,
$F^p(D_t^*,A^*_t)=\dz[[\sigma_{S,V}]\,\big{|}\,\sum_{v\in V}
|a(S,v)|\geq p]$.  We see that $F^p(D_t^*,A^*_t)/F^{p+1}(D_t^*,A^*_t)=
\dz[[\sigma_{S,v}]\,\big{|}\,\sum_{v\in V}|a(S,v)|=p]$, hence
$$E_1^{p,q}(D_t^*,A^*_t)=H^{p+q}(F^p(D_t^*,A^*_t)/F^{p+1}(D_t^*,A^*_t))=
\bigoplus_{a_1,\dots,a_t} H^{p+q}(M_{a_1,\dots,a_t}^*),$$
where the sum is taken over all possible $t$-tuples of
arcs $a_1,\dots,a_t$ such that
\begin{enumerate}
\item[(1)] $a_i\cap\widehat{a_j}=\emptyset$, for any $i\neq j$,
$i,j\in[1,t]$;
\item[(2)] $|a_1|+\dots+|a_t|=p$;
\item[(3)] $\bigcup_{v\in V}\widehat{a(S,v)}\neq V(C_{2r+1})$,
\end{enumerate}
and $M_{a_1,\dots,a_t}^*$ is the~cochain subcomplex generated by all
$\{[\sigma_{S,v}]\}$, such that the arcs with vertices in $V$ are
precisely $a_1,\dots,a_t$, i.e., $\{a(S,v)\,|\,v\in
V\}=\{a_1,\dots,a_t\}$.

Restricting the formula~\eqref{eq:my6} to $M_a^*$, we see that $M_a^*$
is isomorphic to the cochain complex $C^*(\da_{2r-p-2};\dz)$.  More
generally, we see that $M_{a_1,\dots,a_t}^*$ is isomorphic to
$C^*(\da_{2r-\ti p};\dz)$, where $\ti p=\left|\bigcup_{v\in
V}\widehat{a(S,v)}\right|$.

As mentioned, $\ti p\leq 2r$, hence $M_{a_1,\dots,a_t}^*$ is acyclic
for any $a_1,\dots,a_t$ satisfying the above conditions. We conclude
that $(D_t^*,A^*_t)$ is acyclic. The long exact sequence for the
relative cohomology implies that $H^*(A^*_t)=H^*(D^*_t)$. \qed \vspace{5pt}
 
\begin{figure}[hbt]
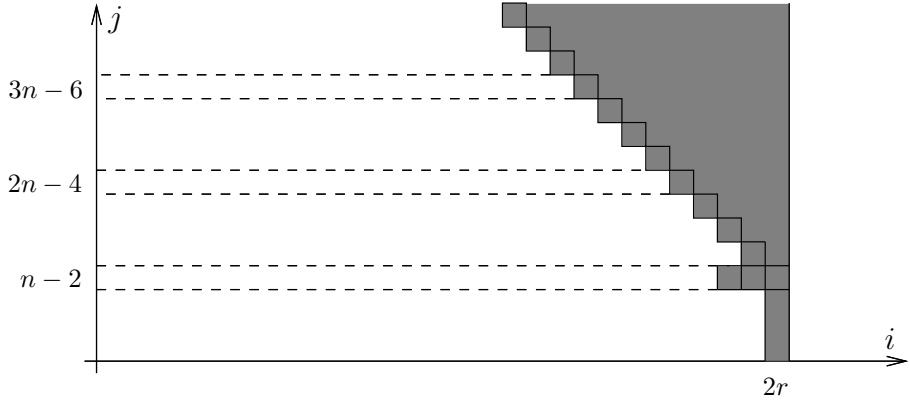

\begin{center}
  \begin{picture}(0,0)%
    \includegraphics{3e1.pstex}%
  \end{picture}%
  \input{3e1.pstex_t}%
 
\end{center}
\caption{The possibly nonzero entries in $E_2^{*,*}$-tableau,
 for $E_2^{p,q}\Rightarrow H^{p+q}(\thomp(C_{2r+1},K_n);\dz)$.}
\label{fig:3e1}
\end{figure}

Now, we can show that $E^{i,t(n-2)}_2=0$ for $t\geq 2$,
$i<2r-(t-1)(n-2)$, that is $E^{*,*}_2$ is $0$ in the region strictly
above row $n-2$ and strictly below the diagonal $x+y=2r+n-2$, see
Figure~\ref{fig:3e1}. Indeed, this is immediate when $2r-2t\geq
2r-(t-1)(n-2)$, which after cancellations reduces to $(n-4)(t-1)\geq
2$. The only cases when this inequality is false are $(t,n)=(2,5)$,
and $n=4$.


\subsection{The case $n=5$, for $E_1^{p,q}\Rightarrow H^{p+q}(\thomp(C_{2r+1},K_n);\dz)$.}
$\,$ \vspace{5pt}

Assume now that $n=5$, $t=2$.

\begin{lm} 
We have $E^{2r-4,6}_2=0$. 
\end{lm}

\pr By dimensional argument, this is true if $2r+1<8$, so we can
assume that $r\geq 4$.

By our previous arguments we need to see that $d_1:A^{2r-4,6}_1\ra
A^{2r-3,6}_1$ is an~injective map. The generators of $A^{2r-4,6}_1$
can be indexed with unordered pairs $\{v,w\}$, $v,w\in V(C_{2r+1})$, such
that 
\[[v-1,v+2]_{2r+1}\cap[w-1,w+2]_{2r+1}=\emptyset,\] 
whereas the generators of $A^{2r-3,6}_1$ can be indexed with ordered
pairs $(v,w)$, $v,w\in V(C_{2r+1})$, such that
\[[v-1,v+2]_{2r+1}\cap[w-1,w+1]_{2r+1}=\emptyset.\] 

In these notations we have
\begin{equation}
\label{eq:n5diff}
d_1(\{v,w\})=\epsilon_1 (v,w)+\epsilon_2 (v,[w+1]_{2r+1})+
\epsilon_3 (w,v)+\epsilon_4 (w,[v+1]_{2r+1}),
\end{equation}
where $\epsilon_1,\epsilon_2,\epsilon_3,\epsilon_4\in\{-1,1\}$.

Take $0\neq\sum_{v,w}\alpha_{v,w}\{v,w\}\in\ker d_1$. Choose $v,w$
such that $\alpha_{v,w}\neq 0$, and the minimum of the two distances
between the arcs $\{v,[v+1]_{2r+1}\}$ and $\{w,[w+1]_{2r+1}\}$ is
minimized. By symmetry we may assume $[w-v-1]_{2r+1}\leq
[v-w-1]_{2r+1}$. Then, it follows from \eqref{eq:n5diff} that
$\alpha_{[v+1]_{2r+1},w}\neq 0$ as well. 

Either $\{[v+1]_{2r+1},w\}$ is not a~well-defined pair or the minimal
distance between the two arcs is smaller for this pair, than for
$\{v,w\}$: $[w-v-1]_{2r+1}\geq [w-v-2]_{2r+1}$. Both ways we get
a~contradiction to the assumption that there exists $\{v,w\}$, such
that $\alpha_{v,w}\neq 0$. We conclude that $d_1:A^{2r-4,6}_1\ra
A^{2r-3,6}_1$ is injective, hence $E^{2r-4,6}_2=0$. \qed \vspace{5pt}

This shows, that when $n\geq 5$, there are no higher differentials
$d_i$, $i\geq 2$, in our spectral sequence, originating in the region
above row $n-2$ and below diagonal $x+y=2r+n-2$. Hence, to figure out
what happens to the entries $E_{\infty}^{2r,n-2}$ and
$E_{\infty}^{2r,n-3}$, it is sufficient to consider rows $n-2$ and
$n-3$.


\begin{figure}[hbt]
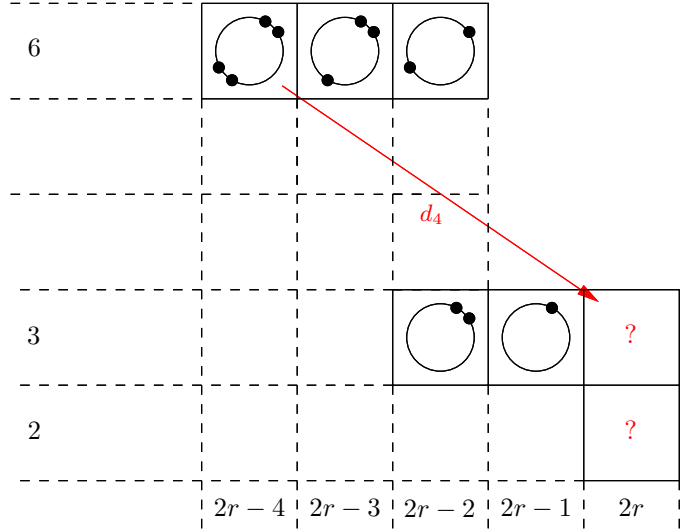

\begin{center}
  \begin{picture}(0,0)%
    \includegraphics{e1n5.pstex}%
  \end{picture}%
  \input{e1n5.pstex_t}%
 
\end{center}
\caption{A part of the $E_1^{*,*}$-tableau, for $n=5$, and 
$E_1^{p,q}\Rightarrow H^{p+q}(\thomp(C_{2r+1},K_n);\dz)$.}
\label{fig:e1n5}
\end{figure}

\subsection{The case $n=4$, for $E_1^{p,q}\Rightarrow H^{p+q}(\thomp(C_{2r+1},K_n);\dz)$.}
$\,$ \vspace{5pt}

For $n=4$ the nonzero rows of $E_1^{*,*}$ are too close to each other,
so we are to do the computation by hand in a~somewhat detailed
way. Since the reduction from $(D_t^*,d_1)$ to $(A_t^*,d_1)$ described
in subsection~\ref{ss4.5} was valid when $n=4$, we may concentrate on
the study of the latter complex. Let us first deal with $(A_2^*,d_1)$. 

\begin{lm}
We have $H^{2r-2}(A_2^*)=H^{2r-3}(A_2^*)=\dz$, and
$H^{i}(A_2^*)=0$, for $i\neq 2r-2,2r-3$.
\end{lm}

\pr We filter $A_2^*$ by
\[F^p A_2^*=\dz\left[\left. [\sigma_{S,\{v_1,v_2\}}]\,\right|\, 
\min(|a(S,v_1)|,|a(S,v_2)|)\geq p\right].\]
Clearly, $\dots\subseteq F^p\subseteq F^{p+1}\subseteq \dots$.
Inspecting the case $p\leq r-2$, we see that in this situation
\[C^*(F^p A_2^*/F^{p-1} A_2^*)=\bigoplus_{i=1}^{2r+1}B_i,\]
where each $B_i$ is isomorphic to $C^*(\Delta_1)$, hence is acyclic.

It follows that $H^*(A_2^*)=H^*(F^{r-1} A_2^*/F^{r-2} A_2^*)$.  Let
$\sigma_i=\sigma_{S_i,V_i}$, where
$S_i=[i+1,i+r-1]_{2r+1}\cup[i+r+2,i-1]_{2r+1}$, $V_i=\{i+1,i+r+2\}$,
and let $\tau_i=\sigma_{\wti S_i,V_i}$, where $\wti
S_i=S_i\cup\{[i+r]_{2r+1}\}$. Clearly,
$d_1(\sigma_i)=\pm\tau_i\pm\tau_{i+r}$, and to verify the statement of
the lemma we need to show that the number of those $\sigma_i$, for
which $d_1(\sigma_i)=\pm(\tau_i+\tau_{i+r})$ is even. By
\eqref{eq:my7} and~\eqref{eq:my6} we see that 
$d_1(\sigma_i)=\pm(\tau_i-\tau_{i+r})$ if $i+r\neq 2r,2r+1$, i.e., the
only cases we need to consider are $i=r$ and $i=r+1$. 

If $i=2r+1$, the different sign comes from \eqref{eq:my6}, and the
sign contribution is $2r+2$. This is an even number, hence again
$d_1(\sigma_i)=\pm(\tau_i-\tau_{i+r})$. If $i=2r$, the different sign
comes from \eqref{eq:my7}, but since $n=4$ is even, the sign remains
the same. \qed \vspace{5pt}

Next, we consider $(A_t^*,d_1)$, for $t\geq 3$.  First, we introduce
some additional notations. Since the sign will not matter in our
argument, we write $\sigma_S$ instead of $\sigma_{S,V}$, it is then
defined only up to a~sign. For $S\subset V(C_{2r+1})$, $\bar
S=V(C_{2r+1})\sm S$; the connected components of $C_{2r+1}[\bar S]$
are called {\it gaps}. Each gap consists of either one or two
elements, we call the first ones {\it singletons}, and the second ones
{\it double gaps}.  Let $m(S)$ be the leftmost element of the gap
which contains $\min(\bar S \cap[2,2r]_{2r+1})$. For $s\in\bar S$, let
$\overleftarrow{s}$ be the leftmost element of the first gap to the
left of the gap containing $s$, and let $\overrightarrow{s}$ be the
leftmost element of the first gap to the right of the gap containing
$s$. For $x,y\in V(C_{2r+1})$, let $d(x,y)$ denote
$\left|[x,y]_{2r+1}\right|-1$.

\begin{lm} \label{lm:4.11}
We have $E_2^{2r-2t,2t}=E_2^{2r-2t+1,2t}=0$.
\end{lm}

\begin{figure}[hbt]
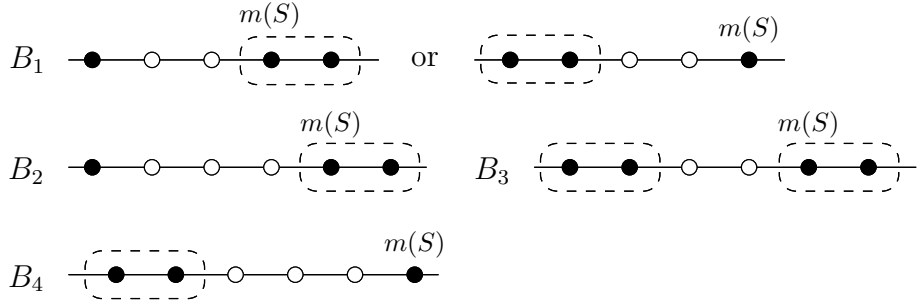

\begin{center}
  \begin{picture}(0,0)%
    \includegraphics{bcases.pstex}%
  \end{picture}%
  \input{bcases.pstex_t}%
 
\end{center}
\caption{The 4 cases in the proof of Lemma~\ref{lm:4.11}.}
\label{fig:bcases}
\end{figure}

\nin{\bf Proof.}\footnote{It is possible to rephrase this argument 
in terms of matchings on chain complexes, see~\cite{dmt}.}  Clearly
$A_t^i=0$, unless $2r-2t\leq i\leq 2r-t$. Note that if $\sigma_S$ is
a~generator of $A_t^{2r-2t+1}$, then $S$ has exactly one singleton. We
decompose $A_t^{2r-2t+1}=B_1\oplus B_2\oplus B_3\oplus B_4$, where
each $B_i$ is spanned by the generators $\sigma_S$, for which certain
conditions are satisfied, see Figure~\ref{fig:bcases}:
\begin{enumerate}
\item[($B_1$)] $\overleftarrow{m(S)}$ is the singleton and
$d(\overleftarrow{m(S)},m(S))=3$, or $m(S)$ is the singleton and
$d(\overleftarrow{m(S)},m(S))=4$;
\item[($B_2$)] $\overleftarrow{m(S)}$ is
the singleton, and $d(\overleftarrow{m(S)},m(S))\geq 4$;
\item[($B_3$)] $m(S)$ and $\overleftarrow{m(S)}$ are in double gaps;
\item[($B_4$)] $m(S)$ is the singleton, and
$d(\overleftarrow{m(S)},m(S))\geq 5$.
\end{enumerate}

Let $(\wti A_t^*,d_1)$ be the complex spanned by $B_1$, $B_2$, $B_3$,
and $A^i_t$, for $2r-2t+2\leq i\leq 2r-t$. The relative complex
$(A_t^*/\wti A_t^*,d_1)$ has only cochains in dimensions $2r-2t$ and
$2r-2t+1$. It is easy to see that the projection $d_1:A_t^{2r-2t}\ra
A_t^{2r-2t+1}/\wti A_t^{2r-2t+1}=B_4$ is an isomorphism, with the
inverse given by $\sigma_S\mapsto \sigma_{S\sm\{m(S)-1\}}$, where
$\sigma_S$ is a~generator from $B_4$. Hence $(A_t^*/\wti A_t^*,d_1)$
is acyclic, and we are led to study the complex $(\wti A_t^*,d_1)$.

Next, we show that $H^{2r-2t+1}(\wti A_t^*)=0$, which is the same as
to say that $d_1$ is injective on $B=B_1\oplus B_2\oplus B_3$.  Let
$\sigma\neq 0$ be in $\ker d_1(B)$. We think of $\sigma$ as a~linear
combination of the generators from the descriptions of $B_1$, $B_2$,
and $B_3$, and let $M$ be the set of generators which have a~nonzero
coefficient in~$\sigma$. 

Assume $M$ contains a generator $\sigma_S$ from $B_3$, and choose
$\sigma_S$ so that $d(m(S),\overrightarrow{m(S)})$ is minimized.  The
coboundary $d_1(\sigma_S)$ contains a~copy of
$\sigma_{S\cup\{m(S)\}}$.  At most~3 other generators will contain
$\sigma_{S\cup\{m(S)\}}$ in the coboundary, depending on which element
we remove from $S\cup\{m(S)\}$ instead of $m(S)$. Since we have chosen
$d(m(S),\overrightarrow{m(S)})$ to be minimal, we cannot remove
$m(S)+2$. Hence, we must remove an~element extending the singleton gap
which is not $m(S)+1$. This gives a~generator of $B_4$, yielding
a~contradiction. 

From now on we may presume that $M$ contains no generators from~$B_4$
or~$B_3$.  Assume now $M$ contains a generator $\sigma_S$ from
$B_2$. Again, choose $\sigma_S$ so that
$d(m(S),\overrightarrow{m(S)})$ is minimized, and note that
$d_1(\sigma_S)$ contains a~copy of $\sigma_{S\cup\{m(S)\}}$.
Examining the generators which contain $\sigma_{S\cup\{m(S)\}}$ in the
coboundary, we see again that, since removing $m(S)+2$ from
$S\cup\{m(S)\}$ would contradict the minimality, we must remove
an~element extending the singleton gap $\overleftarrow{m(S)}$. This
way we will produce a~generator of $B_4$, except for one case: when
$\overleftarrow{m(S)}=1$, and we remove vertex~2. In this case we
produce a~generator from~$B_3$, hence again a~contradiction.

Finally, assume $M$ consists only of generators from $B_1$. Choose
$\sigma_S$ so that $m(S)$ is maximized. Assume first that $m(S)$ is in
a~double gap, and consider the copy of $\sigma_{S\cup\{m(S)\}}$ in
$d_1(\sigma_S)$. There are 3 possibilities. Removing $m(S)+2$ gives
a~generator of $B_2$, whereas removing $m(S)-4$ gives a generator
of~$B_4$. Removing $m(S)-2$ gives either a~generator of $B_3$, if
$m(S)=4$, or a~generator $\sigma_T$ of $B_1$, such that $m(T)=m(S)+1$.
In either case we get a~contradiction.

Assume now that $m(S)$ is a~singleton, and examine the copy of
$\sigma_{S\cup\{\overleftarrow{m(S)}\}}$ in $d_1(\sigma_S)$.  There is
only one possibility for deletion: remove $m(S)+1$.  This will produce
a generator $\sigma_T$ of $B_1$, with $m(T)=m(S)$, but such that
$m(T)$ is in a~double gap, the case which we have already dealt with.

This finishes the proof that $H^{2r-2t+1}(\wti A_t^*)=0$, which,
combined with the acyclicity of $(A_t^*/\wti A_t^*,d_1)$, and the fact
that $H^*(A^*_t)=H^*(D^*_t)$, yields
$E_2^{2r-2t,2t}=E_2^{2r-2t+1,2t}=0$. \qed

\subsection{Finishing the computation of $H^{n-2}(\thom(C_{2r+1},K_n);R)$
and of $H^{n-3}(\thom(C_{2r+1},K_n);R)$,  for
$R=\zz$ or $\dz$.} $\,$ \vspace{5pt} \label{ssect:4.8}

Let us now turn our attention to the cochain complex $A^*_1$. The
generators of $\wti E_1^{2r-2,n-2}$ correspond to arcs of length
$2r-1$ and can be indexed with the elements of $V(C_{2r+1})$: we set
$\tau_{v,2}:=\sigma_{V(C_{2r+1})\sm\{v-2,v-1\},v}$, for any $v\in
V(C_{2r+1})$.  The same way, the generators of $E_1^{2r-1,n-2}$
correspond to arcs of length $2r$, we denote them by setting
$\tau_{v,1}:=\sigma_{V(C_{2r+1})\sm\{v-1\},v}$, for any $v\in
V(C_{2r+1})$.  It follows from~\eqref{eq:my6} that
$$d_1([\tau_{v,2}])=(-1)^{v+1}[\tau_{v,1}]+(-1)^v[\tau_{v-1,1}],$$ for
$v=3,\dots,2r+1$, where the second sign follows from~\eqref{eq:my7};
$$d_1([\tau_{2,2}])=(-1)^{n+1}[\tau_{2,1}]-[\tau_{1,1}],$$
where
the~first sign is determined by the fact that there are $n+2r-3$
vertices before the one inserted at position $2r+1$;
$$d_1([\tau_{1,2}])=(-1)^{n+1}[\tau_{1,1}]+(-1)^n[\tau_{2r+1,1}],$$
where we use again that there are $n+2r-3$ vertices before the
inserted one, and, for determining the second sign, we use the fact
that positions $1$ and $2r+1$ have different parity in
$[1,2r+1]\sm\{2r\}$.

Summarizing, we have the following matrix for the first differential
in $A^*_1$:
$$M=
\begin{bmatrix}
(-1)^{n+1} &          0 & 0 &  0 & \dots & 0 & (-1)^n \\
        -1 & (-1)^{n+1} & 0 &  0 & \dots & 0 & 0 \\
         0 &         -1 & 1 &  0 & \dots & 0 & 0 \\
         0 &          0 & 1 & -1 & \dots & 0 & 0 \\
       \dots & \dots & \dots & \dots & \dots & \dots & \dots \\
         0 &          0 & 0 &  0& \dots & -1 & 1
\end{bmatrix}$$

Assume first that $n\geq 5$, and $(n,r)\neq(5,3)$.

\nin{\bf Case 1:} {\it $n$ is odd.}  It is easy to see that the 
kernel of the differential $d_1:\wti E_1^{2r-2,n-2}\ra E_1^{2r-1,n-2}$
is one-dimensional and is spanned by
\[ [\tau_{1,2}]+[\tau_{2,2}]+[\tau_{3,2}]-[\tau_{4,2}]+[\tau_{5,2}]-
[\tau_{6,2}]+\dots+[\tau_{2r+1,2}],\]
while the image is
\[\left\{\left[\sum_{i=1}^{2r+1}c_i\tau_{i,1}\right]\,\left|\,
\sum_{i=1}^{2r+1}c_i=0\right.\right\}.\]

It follows that $E_2^{2r-2,n-2}=\dz$. Recall that, by
Corollary~\ref{crl:h+}, the cohomology groups of
$\thomp(C_{2r+1},K_n)$ vanish in dimension $n+2r-2$ and less. Hence,
since $d_2:E_2^{2r-2,n-2}\ra E_2^{2r,n-3}$ must be an~isomorphism, we
have $E_1^{2r,n-3}=E_2^{2r,n-3}=\dz$. On the other hand, the map
$d_1:E_1^{2r-1,n-2}\ra E_1^{2r,n-2}$ is surjective, and
$E_2^{2r-1,n-2}=E_2^{2r,n-2}=0$, so $E_1^{2r,n-2}=\dz$.

\nin{\bf Case 2:} {\it $n$ is even.} In this case the map
$d_1:\wti E_1^{2r-2,n-2}\ra E_1^{2r-1,n-2}$ is injective. It follows
that $E_2^{2r-2,n-2}=0$, and hence $E_1^{2r,n-3}=E_2^{2r,n-3}=0$. The
image on the other hand is not the whole $E_1^{2r-1,n-2}$, but only
\[\left\{\left[\sum_{i=1}^{2r+1}c_i\tau_{i,1}\right]\,\left|\,
\sum_{i=1}^{2r+1}c_i\equiv 0\,(\text{mod }2)\right\}.\right.\]
The fact that $E_2^{2r-1,n-2}=E_2^{2r,n-2}=0$ and the surjectivity of
the map $d_1:E_1^{2r-1,n-2}\ra E_1^{2r,n-2}$ imply that
$E_1^{2r,n-2}=\zz$. Again, we used that ${\wti H}^i
(\thomp(C_{2r+1},K_n))$ vanish in dimension $n+2r-2$ and less.

If $(n,r)=(5,3)$, then the argument above essentially holds, with the
exception that $d_1:E_1^{2r-1,n-2}\lra E_1^{2r,n-2}$ does not have to
be surjective. Instead, $\im d_1=\dz$ and $E_1^{2r,n-2}/\im
d_1=\dz$. Thus $H^{n-2}(\thom(C_7,K_5);\dz)=E_1^{2r,n-2}=\dz^2$.

\begin{table}[hbt]
\[\begin{array}{l|c|c|c|} 
(n,r)&\,\,\,\,\,\,R\,\,\,\,\,\,& \,\,\,\,\,\,H^{n-2}\,\,\,\,\,\, &
\,\,\, \,\,\,H^{n-3}\,\,\, \,\,\,\\ 
\hline &&&\\[-0.35cm]
2\not{|}\,\,n, n\geq 5, (n,r)\neq (5,3)\,\,\, \,\,\, & \dz &\dz&\dz\\ 
\hline &&&\\[-0.3cm]
(n,r)=(5,3)& \dz  &\dz^2&\dz\\ \hline
2\,\,|\,\,n, n\geq 6, \text{ or } &&&\\
n=4, r\leq 3 & \dz&\zz&0\\ \hline &&&\\[-0.35cm]
n=4, r\geq 4& \dz  &\dz\oplus\zz&0\\ \hline
n\geq 5, (n,r)\neq (5,3),\text{ or } &&&\\
n=4, r\leq 3& \zz &\zz&\zz\\ \hline 
(n,r)=(5,3), \text{ or } &&&\\
n=4, r\geq 4 & \zz  &{\mathbb Z}_2^2&\zz\\ \hline
\end{array}\]
\caption{$\,$}
\label{tab:answer}
\end{table}

Assume, finally, that $n=4$. If $2r+1=5$, then the computations above
hold. If $2r+1\geq 7$, the argument above still shows that the image
of the map $d_1:E_1^{2r-1,2}\ra E_1^{2r,2}$ is $\dz[d_1(\tau_{1,r})]$,
and $2 d_1(\tau_{1,r})=d_1(2\tau_{1,r})=0$. If $2r+1\geq 9$, we can
compute $E_1^{2r,2}$ and $E_1^{2r,1}$ completely, since
$H^i(\thomp(C_{2r+1},K_n))$ vanish in dimension $n+2r-2$ and less. In
this case, the image of the map $d_1:E_1^{2r-1,2}\ra E_1^{2r,2}$ is
$\zz$, and the map $d_2:E_2^{2r-3,4}\ra E_2^{2r,2}$ is
an~isomorphism. Since, as we have shown earlier, $E_2^{2r-3,4}=\dz$,
we conclude that $E_1^{2r,2}=\dz\oplus\zz$, and $E_1^{2r,1}=0$.

It follows that, for all $(n,r)$, $H^i(\thom(C_{2r+1},K_n);R)=0$, if
$i\in [1,n-4]$, and $R=\dz$ or $\zz$\footnote{This has been
strengthened to yield connectivity in \cite{CK2}, later a shorter
proof appeared in~\cite{En05}.}. We summarize our computations of the
next two cohomology groups in Table~1, where
$H^i=H^i(\thom(C_{2r+1},K_n);R)$, and the case $(n,r)=(4,3)$ is
conjectural.

\noindent
{\bf Proof of Theorem~\ref{thm:even_n}.}

\noindent For $n\geq 6$, and $n=4$, $r\leq 3$, this follows from 
the fact that the target group of the map is~$\zz$. For $n=4$, $r\geq
4$, we have shown above, that $2d_1(\tau_{1,r})=0$. By the
construction, $\tau_{1,r}=\sigma_{V(C_{2r+1})\sm\{r-1\},r}$, so
$d_1(\tau_{1,r})=\pm\sigma_{V(C_{2r+1}),r}$. Let $V(K_2)=\{1,2\}$, and
pick a~nontrivial element $\alpha\in H^{n-2}(\thom(K_2,K_n);\dz)$ by
setting $\alpha:=\eta^*$, $\eta(1):=[1,n-1]$,
$\eta(2):=\{n\}$. Clearly,
$\iota^*_{K_n}(\alpha)=\pm\sigma_{V(C_{2r+1}),r}$, where $\iota(1)=r$,
$\iota(2)=r+1$. Thus, we see that $2\cdot\iota^*_{K_n}(\alpha)=\pm 2
d_1(\tau_{1,r})=0$. \qed

\subsection{The $\zz$-action on the cohomology groups of 
  $\thom(C_{2r+1},K_n)$ for odd~$n$.} $\,$ 

Throughout this subsection we assume that $n$ is odd, and that
$(n,r)\neq (5,3)$. We tensor all our groups with $\mathbb C$ to
simplify the representations. We denote by~$\chi_i$ the
one-dimensional representation of $\zz$ given by the multiplication
by~$(-1)^i$.

\begin{lm}\label{charlm1}
We have $E_1^{2r,n-2}=\chi_r$, as a~$\zz$-module.
\end{lm}
\pr Recall that $\sigma_{V(C_{2r+1}),2r+1}:=\sum_{\eta}\eta^*_+$,
where the sum is taken over all $\eta$, such that
$\eta(2r+1)=[1,n-1]$, and $|\eta(i)|=1$, for all
$i=1,\dots,2r$. $\sigma_{V(C_{2r+1}),2r+1}$ is a~representative of the
generator of $E_1^{2r,n-2}$. Clearly, $\{\eta\circ\gamma\}=\{\eta\}$
as a~collection of cells. To orient the cells in the standard way we
need to reverse $\gamma$ as the permutation of $V(C_{2r+1})$. The sign
of this is $(-1)^r$, hence
$\gamma([\sigma_{V(C_{2r+1}),2r+1}])=(-1)^r[\sigma_{V(C_{2r+1}),2r+1}]$. \qed

\begin{lm}\label{charlm2}
We have $E_1^{2r-1,n-2}=r\chi_0+r\chi_1+\chi_{n+r+1}$, as a~$\zz$-module.
\end{lm}
\pr $\tau_{1,1},\dots,\tau_{2r+1,1}$ can be taken as the
representatives of the generators of $E_1^{2r-1,n-2}$. We see first
that
\begin{equation}
  \label{eq:tau1,1}
  \gamma([\tau_{1,1}])=(-1)^{n+r+1}[\tau_{1,1}].
\end{equation}
Indeed, $\gamma([\tau_{1,1}])=\sgn\pi\cdot[\sigma_{[1,2r],2r}]$,
where $\pi$ is the permutation induced by $\gamma$ on the vertices
of each support simplex of $\tau_{1,1}$, i.e.,
$$\pi=(n+2r-2,n+2r-3,\dots,n+1,n,1,\dots,n-1).$$
Since $\pi$ consists of inverting the sequence $(1,\dots,n+2r-2)$,
and then inverting the subsequence $(1,\dots,n-1)$, we see that
$$\sgn\pi=(-1)^{\left\lfloor\frac{n+2r-2}{2}\right\rfloor+
\left\lfloor\frac{n-1}{2}\right\rfloor}=(-1)^{r-1+\left\lfloor\frac{n}{2}\right\rfloor+
\left\lfloor\frac{n-1}{2}\right\rfloor}=(-1)^{r+n},$$
where we used the fact that the sign of inverting a~sequence
$[1,\dots,m]$ is $(-1)^{\left\lfloor\frac{m}{2}\right\rfloor}$, and that, for any
natural number $m$, we have $\left\lfloor\frac{m}{2}\right\rfloor+
\left\lfloor\frac{m-1}{2}\right\rfloor=m-1$. Additionally,
$[\sigma_{[1,2r],2r}]=-[\sigma_{[1,2r],1}]$ by~(\ref{eq:my7}),
hence~\ref{eq:tau1,1} follows.

Next, we shall see that
\begin{equation}
  \label{eq:taui,1}
  \gamma([\tau_{2r+2-i,1}])=(-1)^{r}[\tau_{i+1,1}],
\end{equation}
for $i=1,\dots,2r$. Again,
$\gamma([\tau_{2r+2-i,1}])=\sgn\pi\cdot[\sigma_{V(C_{2r+1})\sm\{i\},i-1}]$,
where $\pi$ consists of inverting the sequence $(1,\dots,n+2r-3)$, and
then inverting some subsequence of length $n-1$ back. It follows that
\[\sgn\pi=(-1)^{\left\lfloor\frac{n+2r-3}{2}\right\rfloor+
\left\lfloor\frac{n-1}{2}\right\rfloor}=(-1)^{r-1+\left\lfloor\frac{n-1}{2}\right\rfloor+
\left\lfloor\frac{n-1}{2}\right\rfloor}=(-1)^{r+1}.\]
On the other hand, by~(\ref{eq:my7})
$[\sigma_{V(C_{2r+1})\sm\{i\},i-1}]=-[\sigma_{V(C_{2r+1})\sm\{i\},i+1}]$,
hence we get~\eqref{eq:taui,1}. The actual sign has no bearing on our
final conclusion.

Since the permutation action of $\zz$ on a 2-dimensional space
decomposes as $\chi_0+\chi_1$, the formulae~(\ref{eq:tau1,1})
and~(\ref{eq:taui,1}) together yield the claim of the lemma.  \qed

\begin{lm}\label{charlm3}
We have $E_1^{2r-2,n-2}=r\chi_0+r\chi_1+\chi_{r+1}$, as a~$\zz$-module.
\end{lm}
\pr $\tau_{1,2},\dots,\tau_{2r+1,2}$ can be taken as the
representatives of the generators of $E_1^{2r-2,n-2}$. We see first
that
\begin{equation}
  \label{eq:taur,2}
  \gamma([\tau_{r+2,2}])=(-1)^{r+1}[\tau_{r+2,2}].
\end{equation}
We have
$\gamma([\tau_{r+2,2}])=\sgn\pi\cdot[\sigma_{V(C_{2r+1})\sm\{r,r+1\},r-1}]$,
where $\pi$ consists of inverting the sequence of length $n+2r-4$, and
then inverting some subsequence of length $n-1$ back. It follows that
$$\sgn\pi=(-1)^{\left\lfloor\frac{n+2r-4}{2}\right\rfloor+
  \left\lfloor\frac{n-1}{2}\right\rfloor}=(-1)^{r-2+\left\lfloor\frac{n}{2}\right\rfloor+
  \left\lfloor\frac{n-1}{2}\right\rfloor}=(-1)^{r+n+1}.$$
Furthermore,
by~(\ref{eq:my7})
$[\sigma_{V(C_{2r+1})\sm\{r,r+1\},r-1}]=(-1)^n[\sigma_{V(C_{2r+1})\sm\{r,r+1\},r+2}]$,
where $(-1)^n$ is composed of $2r-3$ steps changing the sign, and one step
changing the sign by $(-1)^{n+1}$, since $1$ and $2r+1$ have the same
parity in $V(C_{2r+1})\sm\{r,r+1\}$. Summarizing we get~(\ref{eq:taur,2}).

Second we note that
\begin{equation}
  \label{eq:taui,2}
  \gamma([\tau_{2r+2-i,1}])=\pm[\tau_{i+2,1}],
\end{equation}
for $i\in V(C_{2r+1})\sm\{r+2\}$. Indeed, as before we see that
$\gamma([\tau_{2r+2-i,1}])=\pm[\sigma_{V(C_{2r+1})\sm\{i,i+1\},i-1}]=
\pm[\sigma_{V(C_{2r+1})\sm\{i,i+1\},i+2}]$.

Equations~\eqref{eq:taur,2} and~\eqref{eq:taui,2} show that the
$\zz$-representation splits into $\chi_{r+1}$ and the $r$-fold
permutation action, yielding the claim of the lemma. 
\qed

\begin{crl}
The group $\zz$ acts trivially on
$H^{n-2}(\thom(C_{2r+1},K_n);\dz)=\dz$, and as a~multiplication
by~$-1$ on $H^{n-3}(\thom(C_{2r+1},K_n);\dz)=\dz$.
\end{crl}
\pr It follows from Lemmata~\ref{charlm1}, \ref{charlm2}, and
\ref{charlm3} that $E_1^{2r,n-3}=\chi_{r+1}$, as a~$\zz$-module. The
result follows now from the equation~(\ref{eq:act_disc}). \qed

\section{Cohomology groups of $\zz$-quotients of products of spheres.}

\label{ss5.4}

From now on, unless explicitely stated otherwise, we shall only work
with $\zz$-coefficients.

We begin by introducing another piece of terminology: for a~positive
integer~$d$, let {\it $d$-symbols} be elements of the set
$\{*,\infty\}$, where $*$ will denote a~open $d$-cell, and $\infty$
denote a~$0$-cell. We assume throughout this section that $d\geq
2$. For example, $S^d$ is decomposed into $*$ and $\infty$, whereas
a~direct product of $t$ $d$-dimensional spheres decomposes into cells,
indexed by all possible $t$-tuples of $d$-symbols. We let $\dim *=d$,
$\dim\infty=0$, and we set the dimension of a~tuple of $d$-symbols be
the sum of the dimensions of the constituting symbols.

\subsection{Cohomology groups of $\zz$-quotients of products 
of odd number of spheres.} \label{sssodd} $\,$

Let $X$ be a direct product of $2t+1$ $d$-dimensional spheres, and let
$\zz$ act on $X$ be swapping spheres numbered $2i+1$ and $2i$, for
$i\in[1,t]$, and acting on the first sphere by an antipodal map.  We
shall decompose $X/\zz$ into cells, and describe its cohomology
groups.

Clearly, $X/\zz$ is a total space of a fiber bundle over $\rp^d$ with
fiber homeomorphic to a~direct product of $2t$ $d$-dimensional
spheres.  Consider the standard cell decomposition of $\rp^d$ with one
cell in each dimension $i\in [0,d]$.

\begin{prop}\label{pr:odd}
The space $X/\zz$ can be decomposed into cells indexed with $(i,x,y)$,
where $x$ and $y$ are $t$-tuples of $d$-symbols, $0\leq i\leq d$. The
dimension of this cell is $\dim(i,x,y)=i+\dim x+\dim y$.

The coboundary is given by the equation
\begin{equation}\label{eq:oddcb}
d^{i+\dim x+\dim y}
((i,x,y)^*)=(i+1,x,y)^*+(i+1,y,x)^*,
\end{equation}
where the cochains are considered with $\zz$ coefficients.
\end{prop}
\pr Divide $X$ into cells, by taking the product cell structure, where
spheres 2 to $2t+1$ have one $0$-cell and one $d$-cell, whereas the
first sphere is subdivided as a join of $d+1$ $0$-spheres, with $\zz$
acting antipodally on each of these $0$-spheres. The cells can then be
indexed with triples $(i,x,y)_+$ and $(i,x,y)_-$. The coboundary is
given by
\begin{equation}\label{eq:prcb}
d((i,x,y)_+^*)=(i+1,x,y)_+^*+(i+1,x,y)_-^*.
\end{equation}
This cell structure is $\zz$-equivariant, and no cells are preserved
by the involution. This means that it induces a cell structure on
$X/\zz$. Let $(i,x,y)$ denote the orbit $\{(i,x,y)_+,(i,y,x)_-\}$.
After taking the quotient, \eqref{eq:prcb} becomes~\eqref{eq:oddcb}.
\qed \vspace{5pt}

It follows from the Proposition \ref{pr:odd} that the generators 
of $H^*(X/\zz;\zz)$ are indexed with 
\begin{itemize}
\item $(i,x,x)$, for any $0\leq i\leq d$, and a $t$-tuple of
$d$-symbols $x$, here $(i,x,x)^*$ is the cocycle;
\item $(0,x,y)$, for any $t$-tuples of $d$-symbols $x\neq y$, here
$(0,x,y)^*+(0,y,x)^*$ is the cocycle; $(0,x,y)$ and $(0,y,x)$ index
the same generator;
\item $(d,x,y)$, for any $t$-tuples of $d$-symbols $x\neq y$, here
$(d,x,y)^*$ is the cocycle; $(d,x,y)$ and $(d,y,x)$ index the same
generator.
\end{itemize}

 In other words, the cohomology generators are indexed by pairs
$(\langle A\rangle,i)$, where A is a~$2\times t$ array of $d$-symbols,
and $i\in[0,d]$, if $A$ is fixed by $\zz$, while $i\in\{0,d\}$, if $A$
is not fixed by $\zz$. Here $\zz$ acts on the set of all $2\times t$
arrays of $d$-symbols by swapping the two rows, and $\langle -\rangle$
denotes an~orbit of this action.

For future reference, we remark the following property: these
generators behave functorially, under the maps which insert additional
pairs of spheres. More specifically, assume $q\geq t$, and let
$f:[1,t]\hookrightarrow [1,q]$ be an injection. Let $\ti f:
\underbrace{S^d\times\dots\times S^d}_{2q+1}\ra
\underbrace{S^d\times\dots\times S^d}_{2t+1}$ be the following map:
$\ti f$ is identity on the first sphere, it maps isomorphically the
spheres indexed $2i$ and $2i+1$, for $i\in\im f$, to the spheres
indexed by $2f^{-1}(i)$ and $2f^{-1}(i)+1$, and it maps the remaining
spheres to the base point. Then, the induced map on the cohomology
$\ti f^*$ maps the generator $(\langle A\rangle,i)$ to the generator
$(\langle \wti A\rangle,i)$, where $\wti A$ is the $2\times q$ array
obtained from $A$ as follows: the column $f(i)$ in $\wti A$ is equal
to the column~$i$ in $A$, and, for $j\notin\im f$, the column $j$ in
$\wti A$ consists of two~$\infty$'s.

\subsection{Cohomology groups of $\zz$-quotients of products 
of even number of spheres.}\label{ssseven} $\,$

Let $X$ be a direct product of $2t$ $d$-dimensional spheres, and let
$\zz$ act on $X$ be swapping spheres $2i-1$ and $2i$, for $i\in[1,t]$.
A~customary notation for $X/\zz$ is
$SP^2(\underbrace{S^d\times\dots\times S^d}_{t})$.  Again, we shall
decompose $X/\zz$ into cells, and describe its cohomology groups.

\begin{prop}\label{pr:even}
The space $X/\zz$ can be decomposed into cells indexed with two types
of labels:
\begin{enumerate}
\item [Type 1.] the unordered pairs $\{x,y\}$, where $x$ and $y$ are
$t$-tuples of $d$-symbols, $x\neq y$; the dimension is $\dim x+\dim y$;
\item [Type 2.] $(x,x,k)$, where $x$ is a~$t$-tuple of $d$-symbols,
and $0\leq k\leq\dim x$; the dimension is $\dim  x+k$.
\end{enumerate}

With $\zz$ coefficients, the coboundary is equal to 0 for all
generators, except for $(x,x,0)^*$, when $\dim x\geq 1$, in which case
$d^{\dim x}(x,x,0)^*=(x,x,1)^*$. In particular, the generators of
$\wti H^i(X/\zz;\zz)$ are indexed with the same symbols as the cells
in our decomposition, except for $(x,x,0)$ and $(x,x,1)$.
\end{prop}
\pr Start with a usual subdivision of a direct product of $2t$
$d$-spheres, with the cells indexed by pairs $(x,y)$ of $t$-tuples of
$d$-symbols. For $x\neq y$, the set $(x,y)\cup(y,x)/\zz$ is a~cell in
$X/\zz$, which we label with $\{x,y\}$.

To do the same for $x=y$, we need to take a finer subdivision of
$(x,x)$. Let $(x,x,k)^+$, resp.\ $(x,x,k)^-$, be the set of all points
$\bar\alpha\in{\mathbb R}^{2\dim x}$, $\bar\alpha=(\alpha_i)_{i\in
[2\dim x]}$, such that $\alpha_j=\alpha_{j+\dim x}$, for $k+1\leq
j\leq \dim x$, and $\alpha_k>\alpha_{k+\dim x}$, resp.\
$\alpha_k<\alpha_{k+\dim x}$. Obviously, $(x,x,k)^+$ and $(x,x,k)^-$
are cells, which are mapped to each other by the $\zz$-action. These
cells are different for $k\geq 1$, whereas $(x,x,0)^+=(x,x,0)^-$ is
fixed pointwise.

Set $(x,x,0):=(x,x,0)^+$, and $(x,x,k):=(x,x,k)^+\cup(x,x,k)^-/\zz$,
for $k\geq 1$. The statements about the coboundary map and the
indexing of the cohomology generators follow immediately from our
construction. \qed \vspace{5pt}

Rephrasing Proposition~\ref{pr:even} in the language of arrays, the
generators of $H^*(X/\zz;\zz)$ are indexed with $\zz$-orbits $\langle
A\rangle$ of $2\times t$ arrays of $d$-symbols, with an additional
index $2\leq i\leq \dim A/2$, if $A$ is fixed by the $\zz$-action.
Here $\dim A$ is the sum of the dimensions of all entries of~$A$.

Again, we have functoriality in the following sense: if $q\geq t$, and
$f:[1,t]\hookrightarrow [1,q]$ is an~injection, define $\ti f:
\underbrace{S^d\times\dots\times S^d}_{2q}\ra
\underbrace{S^d\times\dots\times S^d}_{2t}$ analogously to the one in
subsection~\ref{sssodd}.  Then $\ti f^*$ maps $\langle A\rangle$,
resp.\ $(\langle A\rangle,i)$, to $\langle \wti A\rangle$, resp.\
$(\langle \wti A\rangle,i)$, where $\wti A$ is a~$2\times q$ array of
$d$-symbols obtained from $A$ by inserting the columns consisting
entirely of $\infty$'s in the places indexed by $[1,q]\sm\im f$.

\section{Spectral sequence for $H^*(\thomp(C_{2r+1},K_n)/\zz;\zz)$}
\label{sect5}

Next, we would like to show Theorem~\ref{thmmain1}(b). We assume
that $\varpi_1^{n-2}(\thom(C_{2r+1},K_n))\neq 0$, and arrive to
a~contradiction by doing computations in a~spectral sequence, which we
now proceed to set~up.
 
\subsection{$\zz$-equivariant cell decomposition of $\thomp(C_{2r+1},K_n)$.}
$\,$ \vspace{5pt}

\begin{figure}[hbt]
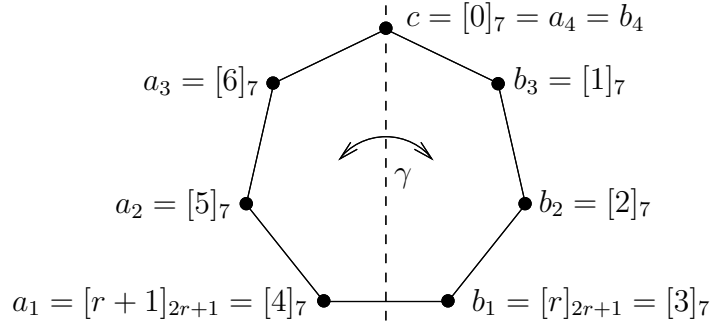

\begin{center}
  \begin{picture}(0,0)%
    \includegraphics{cycle.pstex}%
  \end{picture}%
  \input{cycle.pstex_t}%
 
\end{center}
\caption{Summary of notations.}
\label{fig:cycle}
\end{figure}

For convenience, we give following names to the vertices of
$C_{2r+1}$: $c:=[0]_{2r+1}$, $a_i:=[r+i]_{2r+1}$, and
$b_i:=[r+1-i]_{2r+1}$, for $i\in[1,r]$.  That is $\gamma:C_{2r+1}\ra
C_{2r+1}$ fixes $c$, and $\gamma(a_i)=b_i$, for any $i\in[1,r]$.
Identify $V(C_{2r+1})$ with the vertices of an abstract simplex
$\Delta_{2r}$ of dimension~$2r$. It is also convenient to have
multiple notations for $c$, namely $a_{r+1},b_{r+1}:=c$, see
Figure~\ref{fig:cycle}.

We subdivide the simplex $\Delta_{2r}$ by adding $r$ more vertices,
which we denote $c_1,c_2,\dots,c_r$, and defining a~new abstract
simplicial complex $\tilde\Delta_{2r}$ on the set
$\{c,a_1,\dots,a_r,b_1,\dots,b_r,c_1,\dots,c_r\}=V(\tilde\Delta_{2r})$.
The simplices of $\tilde\Delta_{2r}$ are all the subsets of
$V(\tilde\Delta_{2r})$ which do not contain the subset $\{a_i,b_i\}$,
for any $i\in[1,r]$. We set $\cc=\{c,c_1,\dots,c_r\}$. The complex
$\tilde\Delta_{2r}$ comes equipped with a~simplicial $\zz$-action,
which fixes $\cc$ and swaps $a_i$ and $b_i$, for all $i\in[1,r]$. For
$S\subseteq V(\tilde\Delta_{2r})$ we let $\langle S\rangle$ denote the
$\zz$-orbit of~$S$.


One can think of this new complex $\tilde\Delta_{2r}$ as the one
obtained from $\Delta_{2r}$ by representing it as a~topological join
$\{c\}*[a_1,b_1]*\dots*[a_r,b_r]$, with the additional simplicial
structure defined by inserting an~extra vertex $c_i$ into the middle
of each $[a_i,b_i]$, and then taking the join of $\{c\}$ and the
subdivided intervals.  For $\tilde\sigma\in\tilde\Delta_{2r}$ we
obtain $\vartheta(\tilde\sigma)\in\Delta_{2r}$ by replacing every
$c_i$ in $\ti\sigma$ by $\{a_i,b_i\}$, i.e., $\vartheta(\tilde\sigma)=
(\tilde\sigma\sm\{c_1,\dots,c_r\})\cup\bigcup_{c_i\in\ti\sigma}\{a_i,b_i\}$.

The simplicial complex $\tilde\Delta_{2r}$ has an additional property:
if a~simplex of $\tilde\Delta_{2r}$ is $\gamma$-invariant, then it is
fixed pointwise.  This allows us to introduce a~simplicial structure
(strictly speaking - a~structure of triangulated space) on
$\tilde\Delta_{2r}/\zz$ by taking the orbits of the simplices of
$\tilde\Delta_{2r}$ as the simplices of $\tilde\Delta_{2r}/\zz$.

\vspace{5pt}

\subsection{The chain complex of the subdivision of $\thomp(C_{2r+1},K_n)$.}
$\,$ \vspace{5pt}
 
Since we are working over $\zz$, from now on we shall drop the $+$
notation for the simplices of $\thomp(C_{2r+1},K_n)$, e.g., we shall
write $\eta^*$ instead of $\eta^*_+$ (here we refer to notations
introduced in subsection~\ref{ss3.2}).

Let us now describe a~cochain complex $\wti C^*
(\thomp(C_{2r+1},K_n);\zz)$, which comes from a~triangulation of the
simplicial complex $\thomp(C_{2r+1},K_n)$. The cochain complex
consists of vector spaces over~$\zz$, whose generators are pairs
$(\eta,\sigma)^*$, where $\eta\in\thomp(C_{2r+1},K_n)$, and
$\sigma\in\tilde\Delta_{2r}$, such that
$\vartheta(\sigma)=\supp\eta$. Such a~pair indexes the cochain which
is dual to the cell $\eta\cap\supp^{-1}(\sigma)$. The coboundary of
$(\eta,\sigma)^*$ is the sum of the following generators:
\begin{enumerate}
\item[(1)] $(\ti\eta,\sigma)^*$, if $\supp\ti\eta=\supp\eta$,
$\eta\in\bo\ti\eta$, and $\dim\ti\eta=\dim\eta+1$;
\item[(2)] $(\eta,\sigma\cup\{x\})^*$, if $x\in
V(\ti\Delta_{2r})\sm\sigma$, and $\vt(\sigma)=\vt(\sigma\cup\{x\})$;
\item[(3)] $(\ti\eta,\sigma\cup\{x\})^*$, if $x\in
V(\ti\Delta_{2r})\sm\sigma$, $\ti\eta|_{\vt(\sigma)}=\eta$, and all
the values of $\ti\eta$ on $\vt(\sigma\cup\{x\})\sm\vt(\sigma)$ have
cardinality~1.
\end{enumerate}
 The degree of $(\eta,\sigma)^*$ in $\wti C^*(\thomp(C_{2r+1},K_n))$
is given by
$$\deg(\eta,\sigma)^*=|\sigma|-1+\sum_{v\in\supp\eta}
(|\eta(v)|-1)=\deg\eta+|\sigma|-|\vt(\sigma)|.$$

$\zz$ acts on $\wti C^*(\thomp(C_{2r+1},K_n))$ and we let $\wti
C^*_\zz(\thomp(C_{2r+1},K_n))$ denote its subcomplex consisting of the
invariant cochains. By construction of the subdivision, $\wti
C^*_\zz(\thomp(C_{2r+1},K_n))$ is a~cochain complex for
a~triangulation of the space $\thomp(C_{2r+1},K_n)/\zz$.

\subsection{The filtration of $\wti C^*_\zz(\thomp(C_{2r+1},K_n);\zz)$.}
$\,$ \vspace{5pt}

This time, we consider the natural filtration $(\wti F^{0}
\supseteq\wti F^{1}\supseteq\dots)$ on the cochain complex $\wti
C^*_\zz(\thomp(C_{2r+1},K_n);\zz)$ by the cardinality of~$\sigma$.
Namely, $\wti F^p=\wti F^p C^*_\zz(\thomp(C_{2r+1},K_n);\zz)$, is
a~cochain subcomplex of $C^*_\zz(\thomp(C_{2r+1},K_n);\zz)$ defined
by:
\[\wti F^p: \dots\stackrel{\bo^{q-1}}\lra \wti F^{p,q}\stackrel{\bo^{q}}
\lra\wti F^{p,q+1}\stackrel{\bo^{q+1}}\lra\dots,\] 
where
\[ \wti F^{p,q}=\zz\left[(\eta,\sigma)^*\,\left|\,
(\eta,\sigma)\in C^q_\zz(\thomp(C_{2r+1},K_n);\zz),|\sigma|\geq p+1\right]\right.,
\]
and $\bo^*$ is the restriction of the differential in
$C^*_\zz(\thomp(C_{2r+1},K_n);\zz)$. 

The following formula is the analog of \eqref{eq:e0tab}.

\begin{prop}
For any $p$,
\begin{equation}\label{eq:qe0tab}
\begin{aligned}
\wti F^p/\wti F^{p+1}=&
\bigoplus_{\sigma}C^*(\thom(C_{2r+1}[\vt(\sigma)],K_n)/\zz;\zz)[-p]\\
&\bigoplus_{\langle\tau\rangle}C^*(\thom(C_{2r+1}[\vt(\tau)],K_n);\zz)[-p],
\end{aligned}
\end{equation}
where the first sum is taken over all $\sigma\subseteq\cc$,
$|\sigma|=p+1$, and the second sum is taken over all orbits
$\langle\tau\rangle$, such that $\tau\subseteq V(\ti\Delta_{2r})$,
$|\tau|=p+1$, $\tau\sm\cc\neq\emptyset$.

Hence, the 0th tableau of the spectral sequence associated to
the cochain complex filtration $\wti F^*$ is given by
\begin{equation}
  \label{eq:qE0}
\begin{aligned}
E_0^{p,q}=&
\bigoplus_{\sigma}C^{q}(\thom(C_{2r+1}[\vt(\sigma)],K_n)/\zz;\zz)\\
&\bigoplus_{\langle\tau\rangle}C^{q}(\thom(C_{2r+1}[\vt(\tau)],K_n);\zz),
\end{aligned}
\end{equation}
with the summations over the same sets as in \eqref{eq:qe0tab}.
\end{prop}

\subsection{The analysis of the spectral sequence converging to 
$H^*(\thomp(C_{2r+1},K_n)/\zz;\zz)$.}
$\,$ \vspace{5pt}

The $E_1^{*,*}$-tableau of this spectral sequence is given by
$E_1^{p,q}=H^{p+q}(\wti F^p,\wti F^{p+1})$.  It follows immediately
from the formula \eqref{eq:qE0} that each $E_1^{p,q}$ splits as
a~vector space over $\zz$ into direct sums of
$H^q(\thom(C_{2r+1}[S],K_n);\zz)$, and of
$H^q(\thom(C_{2r+1}[S],K_n)/\zz;\zz)$. More precisely,
\begin{equation}
  \label{eq:qE1}
\begin{aligned}
E_1^{p,q}=&
\bigoplus_{\sigma\subseteq\cc}
H^{q}(\thom(C_{2r+1}[\vt(\sigma)],K_n)/\zz;\zz)\\
&\bigoplus_{\langle\tau\rangle,\,\tau\not\subseteq\cc}
H^{q}(\thom(C_{2r+1}[\vt(\tau)],K_n);\zz).
\end{aligned}
\end{equation}

The generators of $E_1^{p,q}$ stemming from $\sigma\subseteq\cc$ will
be called {\it symmetric}, whereas the generators stemming from
$\langle \tau\rangle$ for $\tau\not\subseteq\cc$ will be called {\it
asymmetric}.

For $i\in[1,r]$, we shall denote the arc
$\{a_i,a_{i-1},\dots,a_1,b_1,b_2,\dots,b_i\}$ by $\smile_i$. For
$2\leq i\leq r$, we denote the arc
$\{a_i,a_{i+1},\dots,a_r,c,b_r,b_{r-1},\dots,b_i\}$ by $\frown_i$.
For $2\leq i<j\leq r$, let $\big{(}_{i,j}$ denote the arc
$\{a_i,a_{i+1},\dots,a_j\}$, let $\big{)}_{i,j}$ denote the arc
$\{b_j,b_{j-1},\dots,b_i\}$, and let $\big{(}\big{)}_{i,j}$ denote
the symmetric pair of arcs $\big{(}_{i,j}$ and~$\big{)}_{i,j}$.

\begin{prop} The map
\begin{equation}\label{eqquot}
q^{n-3}:H^{n-3}(\thom(C_{2r+1},K_n)/\zz;\zz)\ra
H^{n-3}(\thom(C_{2r+1},K_n);\zz),
\end{equation}
is a 0-map.
\end{prop}

\pr First of all, since we are working over the field $\zz$,
the map $q^{n-3}$ is dual to the map on homology
	\[q_{n-3}:H_{n-3}(\thom(C_{2r+1},K_n);\zz)
	\longrightarrow H_{n-3}(\thom(C_{2r+1},K_n)/\zz;\zz),
\]
hence it is enough to prove that $q_{n-3}$ is a $0$-map.

We start by proving that $q_{n-3}=0$ over integers. The map $q_{n-3}$
commutes with the $\zz$-action. Recall that we have proven that
\[H^{n-3}(\thom(C_{2r+1},K_n);\dz)=H^{n-2}(\thom(C_{2r+1},K_n);\dz)=\dz,\] 
so it follows that $H_{n-3}(\thom(C_{2r+1},K_n);\dz)=\dz$. Let $\xi$
be a~generator of the group $H_{n-3}(\thom(C_{2r+1},K_n);\dz)$. By our
previous computations $\gamma^{K_n}(\xi)=-\xi$, since
$H_{n-3}(\thom(C_{2r+1},K_n);\dc)=\chi_1$ as a $\zz$-module (it is
a~dual $\zz$-module to $H^{n-3}(\thom(C_{2r+1},K_n);\dc))$, and since
$H_{n-3}(\thom(C_{2r+1},K_n);\dz)$ is torsion-free.  On the other
hand, the $\zz$-action on $H_{n-3}(\thom(C_{2r+1},K_n)/\zz;\dz)$ is
trivial, hence
\[-q_{n-3}(\xi)=q_{n-3}(-\xi)=q_{n-3}(\gamma^{K_n}(\xi))=
\gamma^{K_n}(q_{n-3}(\xi))=q_{n-3}(\xi).
\]
We conclude that $q_{n-3}(\xi)=0$.

Second, by the universal coefficient theorem the map
\[
\tau:H_{n-3}(\thom(C_{2r+1},K_n);\dz)\otimes\zz\longrightarrow
H_{n-3}(\thom(C_{2r+1},K_n);\zz)
\]
is injective and functorial. In our concrete situation, this map is
also surjective, hence the claim follows from the following
commutative diagram:
$$\begin{CD}
H^{n-3}(\thom(C_{2r+1},K_n);\dz)\otimes\zz   @>\text{0-map}>>
H^{n-3}(\thom(C_{2r+1},K_n)/\zz;\dz)\otimes\zz  \\
        @V\tau V\text{iso}V           @VVV\\
H^{n-3}(\thom(C_{2r+1},K_n);\zz) @>q^{n-3}>>
H^{n-3}(\thom(C_{2r+1},K_n)/\zz;\zz) \\
\end{CD}
$$
$\qquad$\hfill  \qed

\begin{lm}\label{lm:e2r+1n-3}
We have $E_2^{r+1,n-3}=\zz$.
\end{lm}
\pr To start with, the only contribution to $E_1^{r,n-3}$ comes from
$\sigma=\cc$, so the fact that $q^{n-3}$ in \eqref{eqquot} is a~0-map
implies that the differential $d_1:E^{r,n-3}_1\ra E^{r+1,n-3}_1$ is
a~0-map as well.


\begin{figure}[hbt]
\begin{center}
  \begin{picture}(0,0)%
    \includegraphics{2e1.pstex}%
  \end{picture}%
  \input{2e1.pstex_t}%
 
\end{center}
\caption{The $E_2^{*,*}$-tableau, 
$E_2^{p,q}\Rightarrow H^{p+q}(\thomp(C_{2r+1},K_n)/\zz;\zz)$.}
\label{fig:2e1}
\end{figure}

Consider the cochain complex 
\[
A^*:E^{r+1,n-3}_1\stackrel{d_1}{\lra} E^{r+2,n-3}_1\stackrel{d_1}{\lra}
\dots\stackrel{d_1}{\lra} E^{2r,n-3}_1.
\]
The generators of $E_1^{r+i,n-3}$ come from $\tau=\cc\cup I$, for
$I\subseteq\{a_1,\dots,a_r,b_1,\dots,b_r\}$, $|I|=i$. We can identify
the generator indexed by $\langle\tau\rangle$ with the simplex of
${\mathbb R}{\mathbb P}^{r-1}\cong\{a_1,b_1\}*\dots*\{a_r,b_r\}/\zz$,
indexed by $\langle I\rangle$, where the $\zz$-action swaps $a_i$ and
$b_i$, for $i\in[1,r]$.

By inspecting the description of the differential $d_1$ we see that
$A^*$ is isomorphic to the~chain complex $C^*({\mathbb R}{\mathbb
P}^{r-1};\zz)$.  It follows that $E_2^{r+1,n-3}=H^0({\mathbb
R}{\mathbb P}^{r-1};\zz)=\zz$.  \qed \vspace{5pt}

In the proof of the next lemma we shall often use the chain homotopy
between $0$ and the identity.

\noindent
{\it Let $(C^*,d)$ be a cochain complex, and assume there exist
linear maps $\phi^n:C^n\ra C^{n-1}$, $\forall n$, such that
\begin{equation}\label{eq:phid}
\phi^{n+1}(d(\alpha))+d(\phi^n(\alpha))=\alpha, \text{ for all }
\alpha\in C^n.
\end{equation}
Then $C^*$ is acyclic.}

The proof is immediate, since modulo the coboundaries, every
$\alpha\in C^n$ is equal to $\phi^{n+1}(d(\alpha))$, hence
$d(\alpha)=0$ implies $\alpha=0$ modulo the coboundaries.

\begin{lm}\label{lm:e2r-1n-2}
We have $E_2^{r-1,n-2}=0$.
\end{lm}
\pr Set
\[
A^*:E^{0,n-2}_1\stackrel{d_1}{\lra} E^{1,n-2}_1\stackrel{d_1}{\lra}
\dots\stackrel{d_1}{\lra} E^{2r,n-2}_1.
\]
Clearly, to show $E_2^{r-1,n-2}=0$ is the same as to show that
$H^{r-1}(A^*)=0$.

For dimensional reasons, every generator in $A^*$ is indexed either by
$\sigma\subset V(\ti\Delta_{2r})$ with an~arc selected in
$\vt(\sigma)$ (which we call the {\it indexing arc}), or the whole set
$V(C_{2r+1})$ (namely, those coming from
$H^{n-2}(\thom(C_{2r+1},K_n)/\zz)$ and from
$H^{n-2}(\thom(C_{2r+1},K_n))$). To simplify the terminology, we shall
call the set $V(C_{2r+1})$ an~arc as well. Filter the cochain complex
$A^*=G^{2r+1}\supseteq G^{2r}\supseteq\dots\supseteq G^2\supseteq
G^1=0$, where $G^l$ is spanned by the generators whose indexing arc
has the length at least~$l$. We shall compute $H^{r-1}(A^*)$ by
considering the corresponding spectral sequence $\wti E^{p,q}_0:=
C^{p+q}(G^p/G^{p-1})$.

In the same pattern as we have already encountered, the cochain
complex $(G^p/G^{p-1},d_0)$ splits into a~direct sum of subcomplexes
which are indexed by different arcs. For an~arc~$a$, let $B^*_a$
denote the corresponding summand. Hence $\wti E^{p,q}_1= \bigoplus_a
H^{p+q}(B^*_a)$, where the sum is taken over all arcs $a$ of
length~$p$.

Next, by considering all possible arcs case-by-case, we compute the
entries $\wti E^{p,r-1-p}_1$, for $p=2,\dots,2r+1$. To start with,
since $\vt(\sigma)=V(C_{2r+1})$ implies $|\sigma|\geq r+1$, $\wti
E^{2r+1,-r-2}_0=C^{r-1}(G^{2r+1}/G^{2r})=0$ for dimensional reasons,
hence $\wti E^{2r+1,-r-2}_1=0$.

We shall only consider the cases where we cannot use dimensional
reasons to immediately conclude that $B_a^{r-1}=0$.

{\bf Case 1.} Let $a=\smile_r$. Then, $B_a^{r-2}=0$ for dimensional
reasons, and $B_a^{r-1}=\zz$ coming from $\sigma=\{c_1,\dots,c_r\}$.
The differential $d:B_a^{r-1}\ra B_a^r$ is a~0-map since
$f^{n-2}:H^{n-2}(\rp^{n-2};\zz)\ra H^{n-2}(S^{n-2};\zz)$ is a~0-map,
where $f:S^{n-2}\ra S^{n-2}/\zz=\rp^{n-2}$ denotes the covering map.
Hence, in this case, $H^{r-1}(B_a^*)=\zz$.

{\bf Case 2.} Let $a=\frown_2$. Again, $B_a^{r-2}=0$ for dimensional
reasons, and $B_a^{r-1}=\zz$ coming from $\sigma=\{c_2,\dots,c_r,c\}$.
However, this time $d(B_a^{r-1})\neq 0$, since it is induced by the
map $f^{n-2}:H^{n-2}(\thom(G,K_n)/\zz)\ra H^{n-2}(\thom(G,K_n))$,
which, as we have seen, is not a~0-map; here $G$ is the tree on 3
vertices and $\zz$ action is swapping the leaves. Hence
$H^{r-1}(B_a^*)=0$.

{\bf Case 3.} Let $a=\smile_k$, for $1\leq k\leq r-1$.  Let $\alpha\in
B_a^m$ be a~generator indexed by $\sigma\subset V(\ti\Delta_{2r})$. If
$\sigma\subset\cc$, and $x\in\{a_1,\dots,a_r,b_1,\dots,b_r\}\sm
\{a_{k+1},b_{k+1}\}$, then the differential maps $\alpha$ to the
generator indexed by $\langle\sigma\cup\{x\}\rangle$ (again $\smile_k$
is selected) as a~0-map, for the reason described in Case~1. This
means that the complex $B_a^*$ splits into two direct summands, one
containing all generators indexed by $\sigma\subset\cc$, and the other
those indexed by $\langle\sigma\rangle$, such that
$\sigma\sm\cc\neq\emptyset$.

In both summands, define $\phi^m(\alpha)$ to be the generator indexed
by $\langle\sigma\sm\{c\}\rangle$, if $c\in\sigma$, and
$\phi^m(\alpha)=0$ otherwise. The equation \eqref{eq:phid} is
satisfied, so both summands are acyclic, hence so is $B^*_a$, in
particular $H^{r-1}(B_a^*)=0$.

{\bf Case 4.} Let $a=\frown_k$, for $3\leq k\leq r$. We do the same as
in the case~3 with $c$, replaced with $c_{k-2}$. However, in this
complex, if $\sigma\subset\cc$, and $x\in
\{a_1,\dots,a_r,b_1,\dots,b_r\}\sm\{a_{k-1},b_{k-1}\}$, then the
differential maps $\alpha$ to the generator indexed with
$\langle\sigma\cup\{x\}\rangle$ (again $\frown_k$ is selected) as an
identity, hence the complex does not split and the equation
\eqref{eq:phid} can be applied to the whole complex, yielding
$H^{r-1}(B_a^*)=0$.

{\bf Case 5.} Let $a=\big{(}_{2,r}$. For each indexing orbit
$\langle\sigma\rangle$ choose the representative $\sigma$ such that
$a_1\notin\sigma$.  Define $\phi^*$ as in the Case~3, taking $b_1$
instead of $c$. The equation~\eqref{eq:phid} is rather
straightforward. We just need to pay attention to what the
differential does to the generator indexed by
$\sigma=\cc\sm\{c_1,c\}$.

It follows from the description of the cell decomposition and the
cohomology of $SP^2(S^{n-2})$, given as a~special case in
subsection~\ref{ssseven}, that the map
$f^{n-2}:H^{n-2}(SP^2(S^{n-2});\zz)\ra H^{n-2}(S^{n-2}\times
S^{n-2};\zz)$, induced by the quotient map, takes the nonzero element
to the sum of the two generators of $H^{n-2}(S^{n-2}\times
S^{n-2};\zz)$ corresponding to each of the two spheres. In $B_a^*$
this means that the differential of $\sigma$ will contain the
generator of $\langle\sigma\cup\{b_1\}\rangle$, with the indexing arc
$\big{(}_{2,r}$, but not the generator of
$\langle\sigma\cup\{b_1\}\rangle$, with the indexing arc
$\big{)}_{1,r}$. Thus we conclude again that $B^*_a$ is acyclic, and
$H^{r-1}(B_a^*)=0$.

{\bf Case 6.} Let $a=\big{(}_{1,r}$. $B_a^{r-2}=0$ for dimensional
reasons. The space $B_a^{r-1}$ is spanned by the $2^{r-1}$ generators
which we can index with sets $\{a_1,\xi_2,\dots,\xi_r\}$, where
$\xi_i\in\{a_i,c_i\}$, for $2\leq i\leq r$. Denote by $\bar a$ the
generator indexed by the set $\{a_1,a_2,\dots,a_r\}$.  The coboundary
of every generator $\alpha\neq\bar a$ contains some generator $\beta$
indexed by $\{b_i,a_1,\xi_2,\dots,\xi_r\}$. Since $\alpha$ is uniquely
reconstructible from $\beta$, and the coboundary of $\bar a$ does not
contain such generators as $\beta$, we see that an element in
$\ker(d:B_a^{r-1}\ra B_a^r)$ cannot contain $\alpha$ with a~nonzero
coefficient. Thus, the only chance for this kernel to be nontrivial
would be that $\bar a$ lies in it, but, obviously, $d(\bar a)\neq 0$.
Hence, once again, $B^*_a$ is acyclic, and $H^{r-1}(B_a^*)=0$.

{\bf Case 7.} Let $a$ be an~assymetric arc, such that
$a\cap\{a_r,b_r,c\}=\emptyset$. The complex $B^*_a$ is isomorphic to
a~simplicial complex of a cone with an apex in the vertex~$c$. Hence
$B^*_a$ is acyclic, and $H^{r-1}(B_a^*)=0$.

{\bf Case 8.} Let $a$ be such that
$a\cap\{a_1,a_2,b_1,b_2\}=\emptyset$.  The complex $B^*_a$ is
isomorphic to a~simplicial complex of a cone with an apex in the
vertex~$c_1$. Hence $B^*_a$ is acyclic, and $H^{r-1}(B_a^*)=0$.

{\bf Case 9.} Let
$a=\{c,a_i,a_{i+1},\dots,a_r,b_r,b_{r-1},\dots,b_j\}$, for $2\leq
i<j$, where possibly $j=r+1$, which means $a$ does not contain any
$b_i$'s. The complex $B^*_a$ is isomorphic to a~simplicial complex of
a cone with an apex in the vertex~$b_1$. Hence $B^*_a$ is acyclic, and
$H^{r-1}(B_a^*)=0$.

{\bf Case 10.} Let $a=\{a_i,a_{i-1},\dots,a_1,b_1,b_2,\dots,b_j\}$,
for $r\geq i>j\geq 1$. The complex $B^*_a$ is isomorphic to
a~simplicial complex of a cone with an apex in the vertex $a_1$. Hence
$B^*_a$ is acyclic, and $H^{r-1}(B_a^*)=0$.

We can now summarize our computations as follows: $\wti E^{p,r-1-p}_1
=0$, for $p=2,\dots,2r-1$, whereas $\wti E^{2r,-r-1}_1=\zz$.  The
generator of $\wti E^{2r,-r-1}_1$ comes from $a=\smile_r$, which in
turn comes from $\varpi_1^{n-2}(\thom(K_2,K_n))$. The map $d_1:\wti
E^{2r,-r-1}_1\ra \wti E^{2r+1,-r-1}_1$ is the same as 
\[(\iota_{K_n})^{n-2}:H^{n-2}(\thom(K_2,K_n)/\zz;\zz)\ra
H^{n-2}(\thom(C_{2r+1},K_n)/\zz;\zz),\] 
where $\iota:K_2\hookrightarrow C_{2r+1}$ is either of the two
$\zz$-equivariant inclusion maps which take the vertices of $K_2$ to
$\{a_1,b_1\}$.

Since we assumed that $\varpi_1^{n-2}(\thom(C_{2r+1},K_n))\neq 0$, and
the Stiefel-Whitney characteristic classes are functorial, we see that
$d_1:\wti E^{2r,-r-1}_1\ra \wti E^{2r+1,-r-1}_1$ has rank~$1$, hence
$\wti E^{2r,-r-1}_2=0$.  Thus $\wti E^{p,r-1-p}_2 =0$, for
$p=2,\dots,2r+1$, and we conclude that $E_2^{r-1,n-2}=0$.  \qed

\begin{lm}\label{lm:r-i}
We have $E_2^{r-i,n-3+i}=0$, for all $i=2,3,\dots,r$.
\end{lm}
\pr First of all, we note, that the generators in the columns indexed
by $r-1$ and less come from $H^*(\thom(C_{2r+1}[S],K_n)/\zz;\zz)$ and
from $H^*(\thom(C_{2r+1}[S],K_n);\zz)$, with $S\neq V(C_{2r+1})$ in
both cases.

 For each row $q$, $q>n-2$, we shall show that the subcomplex
$A_q^*=(E_1^{*,q},d_1)$ is acyclic in the entry $n+r-q-3$.

We begin by dealing with the case $q=n-1$ separately, that is we
analyze the entry $E_2^{r-2,n-1}$. It follows from
Propositions~\ref{pr:odd}, \ref{pr:even}, and for dimensional reasons,
that the entries $E_1^{0,n-1},E_1^{1,n-1},\dots,E_1^{r-1,n-1}$ are
generated by the contributions whose indexing collections of arcs are
$(\smile_i,\frown_j)$, for $1\leq i<j-1\leq r-1$.

The contributing spaces are homotopy equivalent to $S^{n-2}\times
X/\zz$, where $X$ is a~direct product of $2t+1$ $(n-2)$-dimensional
spheres, and $\zz$-action is as in Section~\ref{ss5.4}. The generators
appearing in the first $r$ entries of the $(n-1)$th row are coming
from the $(n-2)$-cocycle of $S^{n-2}$ and the 1-cocycle of
$\rp^{n-2}$. The analysis of the differentials shows that the complex
$E_1^{0,n-1}\stackrel{d_1}{\ra}E_1^{1,n-1}\stackrel{d_1}{\ra}\dots
\stackrel{d_1}{\ra}E_1^{r-1,n-1}$ computes the nonreduced homology of
a~simplex with $r-2$ vertices (which could be identified with the set
$\{c_2,\dots,c_{r-1}\}$). It follows that the entry $E_2^{r-2,n-1}$,
which computes the first homology group is equal to~$0$.

We assume from now on that $q\geq n$.  Similar to
subsection~\ref{ss4.5} we filter the complexes $A_q^*$.  To describe
the filtration, we sort all generators into 5 groups. The first
group~(Gr1) contains all asymmetric generators, i.e., those coming
from $\langle\sigma\rangle$, for $\sigma\not\subseteq\cc$. The
symmetric generators, coming from $\sigma\subseteq\cc$, are divided
into 4~groups, depending on whether the indexing collection of arcs
\begin{enumerate}
\item[(Gr2)] contains both an $\frown$-arc, and an $\smile$-arc,
\item[(Gr3)] contains an $\smile$-arc, but not an $\frown$-arc,
\item[(Gr4)] contains an $\frown$-arc, but not an $\smile$-arc, 
\item[(Gr5)] contains no $\frown$-arc, and no $\smile$-arc.
\end{enumerate}
The groups are ordered as above. We filter the complex $A_q^*$ by
first sorting the generators by the groups, and then, within each
group we filter additionally by the total length of the indexing arcs.

Let $\wti E^{*,*}_*$ denote the tableaux of the spectral sequence
computing the cohomology of $A_q^*$.  In complete analogy to the
situation in subsection~\ref{ss4.5}, $\wti E^{*,*}_0$ splits into
pieces indexed by various collections of arcs, which we shall call
{\it layers}.

We start by analyzing the contributions of the asymmetric generators.
Consider the subcomplex $B^*$ in the splitting indexed by a~collection
$A$ of $t$ arcs of total length~$l$.  Since the asymmetric generators
come from the direct products of $(n-2)$-spheres, the only nontrivial
cases are $q=t(n-2)$, for $t\geq 2$.

Assume first there is a~gap between some pair of arcs of length at
least~3, and let $x\in V(C_{2r+1})$ be one of the internal points of
a~gap.  If $x=c$, then $B^*$ is isomorphic to the chain complex of
a~cone with apex in~$c$.  Without loss of generality, we can assume
that $x=b_i$, for some $i$.  By the previous assumption,
$b_{i-1},b_{i+1}\notin A$. If $a_i \notin A$, but either $a_{i-1}$, or
$a_{i+1}$ (or both) is in $A$, then $B^*$ is isomorphic to a~chain
complex of a~cone with apex $b_i$. If $a_{i-1},a_i,a_{i+1}\notin A$,
then $B^*$ is isomorphic to a~chain complex of a~cone with apex
$c_i$. Finally, assume $a_i\in A$.  Define $\phi^k:B^k\ra B^{k-1}$ as
follows: for a~generator $\sigma\in B^k$,
\[\phi^k(\sigma)=\begin{cases}
\sigma\setminus\{b_i\}, & \text{ if } b_i\in\sigma,\text{ i.e., if } 
\sigma\cap\{a_i,b_i,c_i\}\text{ is } \{b_i\}, \text{ or }\{b_i,c_i\};\\
\sigma\setminus\{c_i\}, & \text{ if }\sigma\cap\{a_i,b_i,c_i\}=\{a_i,c_i\};\\
0, & \text{ if } \sigma\cap\{a_i,b_i,c_i\}\text{ is } \emptyset, 
\text{ or }\{c_i\},\text{ or }\{a_i\}.
\end{cases}\]
Let $\wti B^*$ be the subcomplex of $B^*$ generated by all $\sigma$,
such that $a_i\in\sigma$. Clearly, the \eqref{eq:phid} is fulfilled
both for $\wti B^*$ and for $B^*/\wti B^*$. It implies that they are
both acyclic, hence so is $B^*$.

If all gaps are of length at most 2, then $l+2t\geq 2r+1$. On the
other hand, $B^p=0$ for $p\leq l/2-1$, since $|\vt(\sigma)|\leq
2|\sigma|-1$, for $\sigma\not\subseteq\cc$. Recall that $q=t(n-2)$, it
follows that the entry $n+r-q-3$ is~0, since
\begin{multline*}
l/2-1-(n+r-t(n-2)-3)>r-t-1-n-r+tn-2t+3=  \\
tn-3t-n+2=(t-1)(n-3)-1\geq 1.
\end{multline*}
Hence $B^*$ is acyclic in the required entry, for all $B^*$ in the
group~(Gr1).

Next, we move on to the symmetric generators. For $\sigma\subseteq\cc$
we call $|\cc\sm\sigma|$ the {\it total length of gaps}. Let $B^*$ be
a~subcomplex in the splitting corresponding to a~layer from the
group~(Gr2). The contributing space here is $S^{n-2}\times X/\zz$,
where $X$ is a~direct product of $2t+1$ $(n-2)$-spheres and
$\zz$-action is as in Section~\ref{ss5.4}.

If $t=0$, since the column number is at most $r-3$, the gap between
the $\smile$-arc and the $\frown$-arc is at least 3. This means that
$B^*$ is isomorphic to a cochain complex of the simplex, hence is
acyclic.

Assume now $t\geq 1$. By examining the cohomology groups of
$S^{n-2}\times X/\zz$, and taking into account that each of the $t$
pairs of spheres must contribute nontrivially, we see that the
dimension of the contributing cocycle of $S^{n-2}\times X/\zz$ is at
least $n-2+t(n-2)=(t+1)(n-2)$, hence the total length of gaps is at
least $t(n-2)+1$. Assume the total length of gaps is at most $2(t+1)$,
as otherwise $B^*$ is isomorphic to a cochain complex of the simplex.
By assumptions, $t\geq 1$ and $n\geq 5$, so unless $(t,n)=(1,5)$, we
have
\[t(n-2)+1-(2t+2)=t(n-4)-1>0,\]
yielding a~contradiction.

Consider the remaining case $(t,n)=(1,5)$. This is the first situation
in which we need to analyze the particular entries of $\wti
E^{*,*}_1$. Since we must have a precise equality, the total length of
gaps is 4, and the only nontrivial case is provided by generators
indexed with $\smile_i$, $\frown_j$, and $\big{(}_{i+3,j-3}$.  The
contributing cohomology generator must be indexed $(0,*,\infty)$, so
just the set of arcs determines everything.

Let $\alpha_{i,j}$ denote such a~generator, and let $\beta_{i,j}$
denote the generator whose indexing set of arcs is $\smile_i$,
$\frown_j$, and $\big{(}_{i+2,j-3}$, and which is also indexed by
$(0,*,\infty)$. Both $\alpha_{i,j}$'s, and $\beta_{i,j}$'s are
generators in $\wti E^{*,*}_1$.  Consider a linear combination
$\sum_{i,j}p_{i,j}\alpha_{i,j}$ lying in the kernel of $d_1$. Since
$d_1(\alpha_{i,j})$ contains $\beta_{i,j}$, $\beta_{i+1,j}$, and no
other $\beta_{i',j'}$'s we see that $p_{i,j}\neq 0$ implies
$p_{i-1,j}\neq 0$.  This leads obviously to $p_{i,j}=0$ for all $i,j$,
hence $B^*$ is acyclic in the required entry.

 Now consider $B^*$ corresponding to a layer from group~(Gr3).  The
contributing space here is $X/\zz$, where $X$ is a~direct product of
$2t+1$ $(n-2)$-spheres and the $\zz$-action is as above. Since we are
in the row $n$ or higher, we must have $t\geq 1$.  The total length of
gaps cannot be larger than $2t+1$, since otherwise $B^*$ is isomorphic
to a~cochain complex of the simplex. On the other hand, since the
dimension of the contributing cohomology generator is at least
$t(n-2)$, the total length of gaps must be at least $(t-1)(n-2)+1$.
Comparing these two we see that
\[(t-1)(n-2)+1-(2t+1)=(t-1)(n-4)-2>0,\]
with exceptions: $t=1$, $n$ is any, $t=2$, $n=5,6$, and $(t,n)=(3,5)$.

Consider first $t=1$. Since we can have at most 3 gaps, we must have
precisely 3 gaps, so the contributing cohomology generators of
$S^{n-2}\times S^{n-2}\times S^{n-2}/\zz$ must have
dimension~$n$. Inspecting the cohomology description of this space
from Section~\ref{ss5.4} we see that there are no generators in
dimensions between $n-2$ and $2n-4$. Since $2n-4>n$ we verify this
case.

Assume now $(t,n)=(2,5)$. The only nontrivial case is when the total
length of gaps is~4 or~5, and $c$ is in the gaps. Let $\alpha_{i,j}$
denote the generator where the gaps are $\{c,i,i+1,j\}$, $r-2\geq
j\geq i+4$, $i\geq 2$, and $\beta_{i,j}$ denote the generator where
the gaps are $\{c,i,j,j+1\}$, $r-3\geq j\geq i+3$, $i\geq 2$. Let
$\gamma_{i,j}$ denote the generator where the gaps are $\{c,i,j\}$,
$r-2\geq j\geq i+3$, $i\geq 2$. Clearly
$d_1(\alpha_{i,j})=\gamma_{i,j}+\gamma_{i+1,j}$, and
$d_1(\beta_{i,j})=\gamma_{i,j}+\gamma_{i,j+1}$. We see that,
restricted to the generators $\alpha_{i,j}$, $\beta_{i,j}$, and
$\gamma_{i,j}$, we have a~chain complex of the graph on Figure~\ref{fig:gr1}.


\begin{figure}[hbt]
\begin{center}
  \begin{picture}(0,0)%
    \includegraphics{gr1.pstex}%
  \end{picture}%
  \input{gr1.pstex_t}%
 
\end{center}
\caption{$\,$}
\label{fig:gr1}
\end{figure}

The kernel is generated by the elementary squares, so it is enough so
see that each square is a coboundary. Indeed, the elementary square
with the lower left corner $(i,j)$ is a~coboundary of the generator
with gaps $\{c,i,i+1,j,j+1\}$.

Finally, assume $(t,n)=(2,6)$ or $(3,5)$. These are the tight cases,
in the sense that the lengths of all gaps are predetermined: the top
gap consists of just $c$, and the other 2, resp.~3, gaps are of
length~2. Assume that the kernel of $d_1$ is not zero, and let
$\alpha$ be an~element in $\ker d_1$. Let $g$ be a~generator, which is
contained in $\alpha$ with a~nonzero coefficient, such that this $g$
maximizes the height of the top gap over all generators appearing with
a~nonzero coefficient in~$\alpha$. Removing the lower element of the
top gap of $g$ gives a~generator which cannot be cancelled out by the
coboundaries of other elements in $\alpha$, due to the assumed
maximality property. This yields a~contradiction, and hence $\ker
d_1=0$.

We move on to group (Gr4), and let $B^*$ correspond to a~generator
indexed by $\frown_j$, and $t$ side arcs.  We can have at most $2t+2$
gaps. The dimension of the contributing cohomology generator is at
least $t(n-2)+n-2$, thus total length of the gaps is at least
$t(n-2)+1$. Comparing these inequalities we get
\[t(n-2)+1-(2t+2)=t(n-4)-1>0,\]
with the only exception $t=1,n=5$.

Let $(t,n)=(1,5)$. The interesting dimension here is 6, thus the total
length of gaps must be exactly~4. The generators $\alpha_i$ indexed
with the collection of arcs $\{\big{(}_{3,i},\frown_{i+3}\}$, for
$4\leq i\leq r-3$.  Since $d_1(\alpha_i)$ contains the generator
indexed with $\{\big{(}_{2,i},\frown_{i+3}\}$, and this generator is
different for different $\alpha_i$, we see that the only linear
combination of $\alpha_i$'s in the kernel of $d_1$ is the trivial one.
Hence, we conclude that the contribution to $\wti E_2^{*,*}$ is~0.

Finally, we consider the case of generators indexed with collections
of arcs avoiding all $\smile$- and $\frown$-arcs. Let us assume there
are $t$ such arcs. To avoid a~cochain complex of a~simplex, the total
length of the gaps must be at most $2t+1$. On the other hand, since
the dimension of the generator is at least $t(n-2)$, the total length
of the gaps must be at least $(t-1)(n-2)+1$. Comparing we see that
\[(t-1)(n-2)+1-(2t+1)=(t-1)(n-4)-2>0,\]
with the exceptions $t=1$, any $n$, $n=5$, $t\leq 3$, and $n=6$,
$t=2$.

If $t=1$, the only nontrivial case occurs in the row $n$. Then, in the
entry of interest we have only one generator: the one indexed by the
arcs $\big{(}\big{)}_{3,r}$. Its coboundary will contain the generator
$\big{(}\big{)}_{2,r}$, hence it is different from~0.

Let $(t,n)=(2,5)$. Since we are in the row 6, for dimensional reasons,
the total length of gaps in the contributing generator is~4. Thus, we
have two types of generators: $\alpha_i^1$ indexed with arc
collections $\{\bbr_{2,i},\bbr_{i+3,r}\}$, $3\leq i\leq r-4$, and
$\alpha_i^2$ indexed with arc collections
$\{\bbr_{3,i},\bbr_{i+2,r}\}$, $4\leq i\leq r-3$. Considering the
value of $d_1$ on the generator indexed with
$\{\bbr_{3,i},\bbr_{i+3,r}\}$, we see that for $i>3$, modulo
coboundaries, any generator $\alpha_i^1$ is a~linear combination of
the generators $\alpha_j^2$. The coboundary of $\alpha_3^1$ contains
the generator indexed with $\{\bbr_{2,3},\bbr_{5,r}\}$, hence no
element in the kernel of $d_1$ can contain $\alpha_3^1$ with a~nonzero
coefficient.  Finally, a~nonzero linear combination of $\alpha_j^2$'s
cannot lie in the kernel of $d_1$, since $d_1(\alpha_j^2)$ contains
the generator indexed with $\{\bbr_{2,j},\bbr_{j+2,r}\}$, which is
different for different~$j$.  Again, we conclude that the contribution
to $\wti E_2^{*,*}$ is~0.

Let $(t,n)=(3,5)$. For dimensional reasons, the total length of the
gaps is precisely~7, thus we have the generators $\alpha_{i,j}$
indexed with arc collections
$\{\bbr_{3,i},\bbr_{i+3,j},\bbr_{j+3,r}\}$, for $4\leq i$, $i+4\leq
j\leq r-4$. Since $d_1(\alpha_{i,j})$ contains the generator indexed
with $\{\bbr_{2,i},\bbr_{i+3,j},\bbr_{j+3,r}\}$, and these generators
are different for different $\alpha_{i,j}$'s, we see that $d_1$ is
injective on the space spanned by $\alpha_{i,j}$'s.  Therefore, in
this case the contribution to $\wti E_2^{*,*}$ is~0. The case
$(t,n)=(2,6)$ is completely analogous.
\qed

\begin{lm} \label{lm6.6}
We have $E_2^{r+i,n-i-1}=0$, for all $i=3,\dots,n-1$.
\end{lm}
\pr Since $|\cc|=r+1$, the entries $E_1^{r+i,n-i-1}$, for
$i=3,\dots,n-1$, come from
$H^{n-i-1}(\thom(C_{2r+1}[\vt(\tau)],K_n))$, for
$\tau\not\subseteq\cc$. We have shown before that these cohomology
groups vanish in dimension $n-4$ and less, which implies
$E_1^{r+i,n-i-1}=0$, hence $E_2^{r+i,n-i-1}=0$.
\qed \vspace{5pt}

We conclude that $E_\infty^{r+1,n-3}=\zz$, contradicting the fact that
$H^{r+n-2}(\thomp(C_{2r+1},K_n)/\zz;\zz)=0$. Therefore, our original
assumption that $\varpi_1^{n-2}(\thom(C_{2r+1},K_n))\neq 0$ is wrong,
and Theorem~\ref{thmmain1}(b) is proved.

\end{document}